\documentclass[10pt, a4paper]{article}
\usepackage[utf8]{inputenc}

\usepackage[T1]{fontenc}
\usepackage{amsmath, amssymb, amsthm}

\theoremstyle{remark}
\newtheorem*{remark}{Remark}
\numberwithin{equation}{subsection}
\usepackage{mathtools}
\usepackage{array}
\usepackage{hyperref}
\usepackage[english]{babel}
\usepackage{booktabs}
\usepackage{latexsym}
\usepackage{bm}
\usepackage[pdftex]{graphicx}
\usepackage{geometry}
\usepackage{breqn}
\usepackage{multicol}
\usepackage[table]{xcolor}
\usepackage{algorithm}
\usepackage{algpseudocode}

\pdfoutput=1
\usepackage{etoolbox}
\AtBeginEnvironment{tabular}{\linespread{0.8}\selectfont}

\newcommand{\mstr}[4]{\prescript{#1}{#2}{!}^{#3}_{#4}{}}
\newcommand{\merge}{\mathbin{@}}

\title{A Formalism for Merging Power Series in Solving Multinomial Equations: Master Series and Master Numbers}
\author{Sergey V. Berezin}
\date{Version 4.0, \today}
\begin{document}

\maketitle

\begin{center}
\textbf{Affiliation}: Municipal Budgetary Institution, Iglino, Republic of Bashkortostan, Russia \\
\textbf{ORCID}: \href{https://orcid.org/0000-0001-8086-8288}{0000-0001-8086-8288} \\
\textbf{Email}: \href{mailto:bkcru@bk.ru}{bkcru@bk.ru} \\
\textbf{License}: Creative Commons Attribution 4.0 International (CC BY 4.0)
\end{center}

\begin{abstract}
It is shown that all solutions of the canonical equation $y^a = 1 + axy^b$ with arbitrary complex exponents $a, b \in \mathbb{C}$ (case $m=1$) and of its limiting cases $y = (1 + ax)^{1/a}$ at $b=0$, $y = e^x$ at $a=0, b=0$, and $y = e^{x y^b}$ at $a=0$, as well as the logarithms of these solutions (case $m=0$), are expressed by a single power master series
\[
M(m; a; b; x) = m + x + \frac{m-a+2b}{2!}\,x^2 + \frac{(m-a+3b)(m-2a+3b)}{3!}\,x^3 + \cdots
\]
For multivalued solutions beyond the radius of convergence, a geometric criterion for selecting the conjugate series (the same master series, but with a different set of parameters $M(m; \tilde{a}; \tilde{b}; \tilde{x})$) is proposed, ensuring branch continuity on the complex plane. This representation of power series in terms of master numbers not only allows one to compute coefficients in degenerate cases where the gamma-function representation loses meaning, but also provides a transparent analytical algorithm for the merge operation on master series. The method is illustrated on problems in gas dynamics and plasma physics, where it is demonstrated that physically realizable states require the application of conjugate series beyond the radius of convergence of the basic expansion.
\end{abstract}


\tableofcontents
\section{Introduction}
\label{sec:intro}
\textbf{Master-Join methods} reduce the solution of trinomial and multinomial equations to a unified formalism of power series based on the master series.
Unlike classical analytic continuation methods (contour integrals, Padé approximants, continuation along chains of disks), which require cumbersome computational procedures, the Master-Join approach transforms the problem to a simple parameter substitution within the same power series: $M(m; a; b; x) \to M(m; \tilde{a}; \tilde{b}; \tilde{x})$.
\subsection*{On the structure of the paper}
In the \textbf{first section}, we consider the canonical trinomial equation
\begin{equation}
y^a = 1 + a x y^b.
\label{eq:intro_main}
\end{equation}
We show how to find its solution beyond the radius of convergence without losing branch continuity.

We introduce the \textbf{master series}
\begin{equation}
M(m; a; b; x) = m + x + \frac{m-a+2b}{2!}x^2 + \frac{(m-a+3b)(m-2a+3b)}{3!}x^3 + \cdots,
\label{eq:intro_master_series}
\end{equation}
which is a single-valued solution of the \textbf{lucky equation}
\begin{equation}
y = (1 + a x y^b)^{1/a} = M(1; a; b; x).
\label{eq:intro_lucky}
\end{equation}

For the general equation
\begin{equation}
p Y^a = q + z Y^b
\label{eq:intro_general}
\end{equation}
we obtain the multivalued solution:
\begin{equation}
Y = e^f M(1; a; b; x), \qquad x = \frac{z e^{b f}}{a q}, \qquad f = \frac{\ln|q/p| + i\bigl(\arg(q/p) + 2\pi n\bigr)}{a}, \qquad n \in \mathbb{Z}.
\label{eq:intro_multivalued}
\end{equation}

The radius of convergence of the master series is given by:
\begin{equation}
R = \frac{|1 - a/b|^{b/a}}{|b - a|}.
\label{eq:intro_radius}
\end{equation}

However, the original general trinomial equation $S_2 Y^{\beta_2} + S_1 Y^{\beta_1} + S_0 = 0$ can be reduced to the canonical form $p Y^a = q + z Y^b$ in three different ways. If the argument of the series exceeds the radius of convergence ($|x| \ge R$) under the base transformation and the series diverges, the other two transformations are employed to find the solution. To uniquely determine which of these provides a continuous analytic continuation of the specified branch $n$, a developed geometric criterion is applied. The essence of the algorithm is extremely simple: it partitions the imaginary axis into sectors corresponding to the branch indices $n$, and selects as the analytic continuation the conjugate series (with index $h$ or $k$) whose sector center falls inside sector $n$. The special function --- the \textbf{real Master ultra-radical} --- is then investigated, formulas for its derivative and integrals are derived, and examples of applications in gas dynamics and plasma physics are presented.

In the \textbf{second section}, the focus shifts to the algebraic structure of the coefficients themselves. We study \textbf{master numbers} --- the fundamental "building blocks" of master series --- deriving their key properties and identities. A central role is played by the introduction of the \textbf{merge} operation for master numbers and series.

In the following sections, it is demonstrated that the Taylor series for the ultraradical in the neighborhood of an arbitrary complex point $x_0$ with displacement $z = x - x_0$ takes the explicit form:
\begin{equation}
y(x_0 + z) = y(x_0) \left[ 1 + \sum_{\ell=1}^{\infty} \frac{P_{\ell-1}(\xi_0; 1)}{\ell!} \psi(x_0)^\ell z^\ell \right],
\label{eq:intro_taylor_factorized}
\end{equation}
where $\psi(x_0) = y'(x_0)/y(x_0)$ is the value of the logarithmic derivative, $\xi_0 = (a-b)b\,x_0\,\psi(x_0)$ is the dimensionless phase-space variable at the expansion point, and the family of derivative polynomials $P_{\ell-1}(\xi_0; 1)$ is computed recursively:
\begin{equation}
P_\ell(\xi; m) = \sum_{k=0}^{\ell} \delta_{\ell,k}(m) \, \xi^k.
\label{eq:intro_P}
\end{equation}
The main methodological result of formula \eqref{eq:intro_taylor_factorized} is the complete elimination of the need to compute complicated higher-order derivatives in explicit form. Instead of cumbersome differentiation, a purely combinatorial three-term recurrence relation --- the \textbf{``Qwen's Triad''} --- is introduced, which computes all polynomial coefficients row by row:
\begin{equation}
\boxed{\delta_{\ell+1,k} = \bigl[m + (\ell+1+k)\sigma\bigr] \delta_{\ell,k} - (\ell+k) \delta_{\ell,k-1} + \tau(k+1) \delta_{\ell,k+1},}
\label{eq:intro_delta}
\end{equation}
where $\sigma = 2b - a$ and $\tau = (a-b)b$ are the parameters determined by the exponents of the original equation.

A fundamental special case of this structure at $x_0 = 0$ (when the phase variable $\xi_0 = 0$) is the reduction of the polynomials to their constant terms, which coincide exactly with the classical master numbers, thereby naturally recovering the base Maclaurin series. Since both expansions describe the same analytic function, all identities and structural properties of the master series are automatically inherited by the Taylor series.

\section{The Master Ultra-Radical and Its Lucky Series}
\label{sec:lucky_equation}

Consider the canonical unit equation
\begin{equation}
y = (1 + axy^b)^{1/a}, \qquad a, b \in \mathbb{C}, \quad x \in \mathbb{C}.
\label{eq:lucky}
\end{equation}
We shall hereafter refer to \eqref{eq:lucky} as the \textbf{lucky equation}.

\subsection{Expanding the lucky equation into a power series}
\label{subsec:derivation}

The power series for the lucky equation can be obtained using the Lagrange inversion formula~\cite{Graham1994, Corless1996, Gessel2016}, for which modern elementary inductive proofs are also available~\cite{SuryaWarnke2023}. However, due to the simplicity of the structure~\eqref{eq:lucky}, it is \emph{advisable} to apply a \emph{combined expansion}: the generalized binomial theorem for complex exponents~\cite{DLMF, Whittaker1927, Graham1994} and the multinomial theorem for integer powers~\cite{Graham1994, Knuth1997, Knuth2002}. This approach also allows one to identify the intermediate convergence conditions at each step of the expansion.

\subsubsection*{Introducing the Kernel}

At $x = 0$, equation \eqref{eq:lucky} yields $y = 1$. Therefore, we seek a solution in the form
\begin{equation}
y = 1 + K, \qquad K = \sum_{\ell=1}^{\infty} \frac{N_\ell}{\ell!}\,x^\ell,
\label{eq:kernel_def}
\end{equation}
where $K$ is the \textbf{kernel} of the master series (a power series without a constant term), and $N_\ell$ are the coefficients to be determined.

\subsubsection*{Step 1: Expansion of the right-hand side}

We apply the generalized binomial theorem to $(1 + axy^b)^{1/a}$. This expansion is valid under the condition:
\begin{equation}
|a x y^b| < 1.
\label{eq:binomial_convergence}
\end{equation}
\textbf{Remark.} For $b=0$, condition \eqref{eq:binomial_convergence} reduces to $|ax| < 1$, which corresponds to the radius of convergence $R = |1/a|$ for the binomial root $(1+ax)^{1/a}$.

\begin{align}
(1 + axy^b)^{1/a}
&= 1 + \frac{1}{a}\,axy^b + \frac{\frac{1}{a}\!\left(\frac{1}{a}-1\right)}{2!}\,a^2x^2y^{2b}
+ \frac{\frac{1}{a}\!\left(\frac{1}{a}-1\right)\!\left(\frac{1}{a}-2\right)}{3!}\,a^3x^3y^{3b} + \cdots
\notag \\
&= 1 + xy^b + \frac{1-a}{2!}\,x^2y^{2b} + \frac{(1-a)(1-2a)}{3!}\,x^3y^{3b}
+ \frac{(1-a)(1-2a)(1-3a)}{4!}\,x^4y^{4b} + \cdots
\label{eq:rhs_binomial}
\end{align}

\subsubsection*{Step 2: Powers of the Kernel}

Substituting $y^{\ell b} = (1+K)^{\ell b}$ into \eqref{eq:rhs_binomial} requires two steps. First, we expand the powers of the kernel $K^j$ using the multinomial theorem:
\begin{align}
K &= \frac{N_1}{1!}x + \frac{N_2}{2!}x^2 + \frac{N_3}{3!}x^3 + \frac{N_4}{4!}x^4 + \cdots,
\label{eq:K1} \\
K^2 &= N_1^2 x^2 + N_1 N_2 x^3 + \left(\frac{N_2^2}{4} + \frac{N_1 N_3}{3}\right)x^4 + \cdots,
\label{eq:K2} \\
K^3 &= N_1^3 x^3 + \frac{3}{2}N_1^2 N_2 x^4 + \cdots,
\label{eq:K3} \\
K^4 &= N_1^4 x^4 + \cdots.
\label{eq:K4}
\end{align}
\subsubsection*{Step 3: Binomial expansion}
Then, we expand $y^{\ell b} = (1+K)^{\ell b}$ using the generalized binomial theorem. This expansion is valid under the condition:
\begin{equation}
|K| < 1 \quad \text{or} \quad b = 0.
\label{eq:kernel_convergence}
\end{equation}
\textbf{Remark.} Condition \eqref{eq:kernel_convergence} highlights the fundamental distinction between the cases $b = 0$ and $b \neq 0$. For $b = 0$, the power series converge for arbitrary $|K|$, however condition \eqref{eq:binomial_convergence} must still hold: $|ax| < 1$, which corresponds to the radius of convergence $R = |1/a|$ for the binomial root $(1+ax)^{1/a}$. For $a = 0$, $b = 0$ (the exponential), the radius of convergence is infinite.
\begin{equation}
y^{\ell b} = 1 + \ell b K + \frac{\ell b(\ell b-1)}{2!}K^2 + \frac{\ell b(\ell b-1)(\ell b-2)}{3!}K^3 + \frac{\ell b(\ell b-1)(\ell b-2)(\ell b-3)}{4!}K^4 + \cdots
\label{eq:y_lb}
\end{equation}

Substituting \eqref{eq:K1}--\eqref{eq:K4} into \eqref{eq:y_lb} and grouping by powers of $x$, we obtain:
\begin{align}
y^{\ell b} &= 1 + \ell b N_1 x
\notag \\
&\quad + x^2\,\frac{\ell b N_2 + \ell b(\ell b-1)N_1^2}{2!}
\notag \\
&\quad + x^3\left(\frac{\ell b N_3}{3!} + \frac{\ell b(\ell b-1)N_1 N_2}{2!} + \frac{\ell b(\ell b-1)(\ell b-2)N_1^3}{3!}\right)
\notag \\
&\quad + x^4\left(\frac{\ell b N_4}{4!} + \frac{\ell b(\ell b-1)}{2!}\!\left(\frac{N_2^2}{4}+\frac{N_1 N_3}{3}\right)
+ \frac{\ell b(\ell b-1)(\ell b-2)N_1^2 N_2}{4}
+ \frac{\ell b(\ell b-1)(\ell b-2)(\ell b-3)N_1^4}{4!}\right) + \cdots
\label{eq:y_lb_expanded}
\end{align}

\subsubsection*{Step 4: Substitution and Grouping}

We substitute the expansions of $y^b$, $y^{2b}$, $y^{3b}$, $y^{4b}$ (from \eqref{eq:y_lb_expanded} for $\ell=1,2,3,4$) into the right-hand side of \eqref{eq:rhs_binomial} and group by powers of $x$:
\begin{align}
(1+axy^b)^{1/a} &= 1 + x
\notag \\
&\quad + x^2\left(bN_1 + \frac{1-a}{2!}\right)
\notag \\
&\quad + x^3\left(\frac{bN_2 + b(b-1)N_1^2}{2!} + (1-a)bN_1 + \frac{(1-a)(1-2a)}{3!}\right)
\notag \\
&\quad + x^4\Bigg(\frac{bN_3}{3!} + \frac{b(b-1)N_1 N_2}{2!} + \frac{b(b-1)(b-2)N_1^3}{3!}
\notag \\
&\qquad\qquad + \frac{(1-a)\bigl(bN_2 + b(2b-1)N_1^2\bigr)}{2}
+ \frac{(1-a)(1-2a)bN_1}{2}
+ \frac{(1-a)(1-2a)(1-3a)}{4!}\Bigg) + \cdots
\label{eq:rhs_grouped}
\end{align}

\subsubsection*{Step 5: Equating Coefficients}

The left-hand side of the lucky equation \eqref{eq:lucky} is $y = 1 + K$, that is,
\begin{equation}
y = 1 + N_1 x + \frac{N_2}{2!}x^2 + \frac{N_3}{3!}x^3 + \frac{N_4}{4!}x^4 + \cdots
\label{eq:lhs}
\end{equation}

Equating the coefficients for equal powers of $x$ in \eqref{eq:lhs} and \eqref{eq:rhs_grouped}, we obtain the system of equations:
\begin{align}
x^1:&\quad N_1 = 1, \label{eq:sys1} \\
x^2:&\quad \frac{N_2}{2!} = bN_1 + \frac{1-a}{2!}, \label{eq:sys2} \\
x^3:&\quad \frac{N_3}{3!} = \frac{bN_2 + b(b-1)N_1^2}{2!} + (1-a)bN_1 + \frac{(1-a)(1-2a)}{3!}, \label{eq:sys3} \\
x^4:&\quad \frac{N_4}{4!} = \frac{bN_3}{3!} + \frac{b(b-1)N_1 N_2}{2!} + \frac{b(b-1)(b-2)N_1^3}{3!}
\notag \\
&\qquad\quad + \frac{(1-a)\bigl(bN_2 + b(2b-1)N_1^2\bigr)}{2}
+ \frac{(1-a)(1-2a)bN_1}{2}
+ \frac{(1-a)(1-2a)(1-3a)}{4!}. \label{eq:sys4}
\end{align}

\subsubsection*{Step 6: Solving the System}

\textbf{Coefficient $N_1$.} From \eqref{eq:sys1}:
\begin{equation}
\boxed{N_1 = 1.}
\end{equation}

\textbf{Coefficient $N_2$.} From \eqref{eq:sys2} with $N_1 = 1$:
\begin{equation}
N_2 = 2b + (1-a) \quad\Longrightarrow\quad \boxed{N_2 = 1 - a + 2b.}
\end{equation}

\textbf{Coefficient $N_3$.} Multiplying \eqref{eq:sys3} by $3! = 6$ and substituting $N_1 = 1$, $N_2 = 1-a+2b$:
\begin{align}
N_3 &= 3\bigl(bN_2 + b(b-1)\bigr) + 6(1-a)b + (1-a)(1-2a)
\notag \\
&= (1-a+3b)(1-2a+3b).
\end{align}
\begin{equation}
\boxed{N_3 = (1-a+3b)(1-2a+3b).}
\end{equation}

\textbf{Coefficient $N_4$.} Multiplying \eqref{eq:sys4} by $4! = 24$ and substituting the found $N_1, N_2, N_3$:
\begin{align}
N_4 &= 4bN_3 + 12b(b-1)N_2 + 4b(b-1)(b-2)
\notag \\
&\quad + 12(1-a)\bigl(bN_2 + b(2b-1)\bigr) + 12(1-a)(1-2a)b + (1-a)(1-2a)(1-3a)
\notag \\
&= (1-a+4b)(1-2a+4b)(1-3a+4b).
\end{align}
\begin{equation}
\boxed{N_4 = (1-a+4b)(1-2a+4b)(1-3a+4b).}
\end{equation}

\subsubsection*{Step 7: General formula}

Analyzing the structure of the obtained coefficients, we notice a pattern: each $N_\ell$ is a product of $\ell-1$ linear factors of the form $(1 - a\gamma + b\ell)$, where $\gamma = 1, 2, \ldots, \ell-1$.
\begin{equation}
N_\ell = \prod_{\gamma=1}^{\ell-1}(1 - a\gamma + b\ell), \qquad \ell \geq 2, \quad N_0 = N_1 = 1.
\label{eq:master_numbers_m1}
\end{equation}

Thus, the lucky equation \eqref{eq:lucky} has a solution in the form of the \textbf{lucky series}, convergent for $|x| < R$:
\begin{equation}
y = (1 + axy^b)^{1/a} = 1 + \sum_{\ell=1}^{\infty} \frac{N_\ell}{\ell!} x^\ell, \qquad |x| < R,
\label{eq:master_series_lucky}
\end{equation}
where the radius of convergence $R$ will be determined below.

In this formula, the structure of each coefficient is clearly visible. The only element that remains outside this structure so far is the constant term $1$. Let us denote it by $m$ and consider the case $m=0$.
\begin{equation*}
0 + \sum_{\ell=1}^{\infty} \frac{x^\ell}{\ell!} \prod_{\gamma=1}^{\ell-1} (0 - a\gamma + b\ell).
\end{equation*}

\subsection{Logarithmic identity: the influence of parameter m on the solution}
\label{subsec:log_identity}

Now consider the same power series, but with parameter $m=0$:
\begin{equation}
w = x + \frac{2b-a}{2!}x^2 + \frac{(3b-a)(3b-2a)}{3!}x^3 + \frac{(4b-a)(4b-2a)(4b-3a)}{4!}x^4 + \cdots
\label{eq:series_m0_}
\end{equation}

We shall show that the exponential of this series yields the lucky series \eqref{eq:master_series_lucky}. Substituting $w$ into the expansion of the exponential $e^w = 1 + w + \frac{w^2}{2!} + \frac{w^3}{3!} + \cdots$:

\textbf{Term at $x^1$:}
\begin{equation}
[w]_{x^1} = x \quad \Rightarrow \quad 1 \cdot x
\end{equation}

\textbf{Term at $x^2$:}
\begin{equation}
[w]_{x^2} + \frac{1}{2}[w^2]_{x^2} = \frac{2b-a}{2}x^2 + \frac{1}{2}x^2 = \frac{1-a+2b}{2}x^2
\end{equation}

\textbf{Term at $x^3$:}
\begin{align}
[w]_{x^3} + \frac{1}{2}[w^2]_{x^3} + \frac{1}{6}[w^3]_{x^3} 
&= \frac{(3b-a)(3b-2a)}{6}x^3 + \frac{1}{2}(2b-a)x^3 + \frac{1}{6}x^3 \notag \\
&= \frac{9b^2 - 9ab + 2a^2 + 6b - 3a + 1}{6}x^3 \notag \\
&= \frac{(1-a+3b)(1-2a+3b)}{6}x^3
\end{align}

Thus, we have obtained exactly the lucky series \eqref{eq:master_series_lucky}:
\begin{equation}
e^{M(0;a;b;x)} = M(1;a;b;x).
\label{eq:exp_identity}
\end{equation}

This \textbf{fundamental exponential identity} shows that series with different values of the parameter $m$ are related through the exponential function.

Here $M(m;a;b;x)$ is the \textbf{master series} --- a universal power series:
\begin{equation}
M(m;a;b;x) = m + \sum_{\ell=1}^{\infty} \frac{x^\ell}{\ell!} \prod_{\gamma=1}^{\ell-1} (m - a\gamma + b\ell).
\label{eq:master_series_def}
\end{equation}

We shall call a \textbf{master number} the natural generalization of factorials, whose identities and properties will be analyzed in the following sections:
\begin{equation}
N(m; a; b; \ell) = \prod_{\gamma=1}^{\ell-1}(m - a\gamma + b\ell).
\label{eq:master_numbers_}
\end{equation}

\textbf{Historical remark.} Belkić~\cite{Belkic2019} obtained a logarithmic series for the trinomial root (formula (12.3) in his paper), expressing it through Pochhammer symbols. This series is a special case of the logarithmic master series $M(0; a; b; x)$ with $a = 1$. The master-number representation generalizes Belkić's result to arbitrary complex $a, b \in \mathbb{C}$ and reveals its connection with the Lambert $W$ function and the Euler $T$ function.

\subsection{First remarkable property: the limiting case b=0}
\label{subsec:b_zero}

Let us write the lucky equation and its logarithmic form in compact notation:
\begin{align}
y &= (1 + axy^b)^{1/a} = M(1;a;b;x), \label{eq:lucky_power} \\
w &= \frac{\ln(1 + ax e^{bw})}{a} = M(0;a;b;x). \label{eq:lucky_log}
\end{align}

Both equations form a natural pair related by the exponential identity $M(1;a;b;x) = \exp(M(0;a;b;x))$.

Setting $b = 0$ in both equations, we have $y^b = y^0 = 1$ and $e^{bw} = e^0 = 1$, and both equations reduce to a remarkably simple form:
\begin{align}
y &= (1 + ax)^{1/a} = M(1;a;0;x), \qquad |x| < R, \label{eq:binomial_limit} \\
w &= \frac{\ln(1 + ax)}{a} = M(0;a;0;x), \qquad |x| < R. \label{eq:log_limit}
\end{align}

Thus, at $b = 0$, the lucky equation \emph{without any additional transformations} reduces to:
\begin{itemize}
\item The \textbf{binomial root} $(1 + ax)^{1/a}$ --- for the power-type case $m=1$;
\item The \textbf{generalized logarithm} $\frac{1}{a}\ln(1 + ax)$ --- for the logarithmic case $m=0$.
\end{itemize}

\subsection{Second remarkable property: the limiting case a=0}
\label{subsec:a_zero}
Consider the limiting transition $a \to 0$ in the lucky equation and its logarithmic form.

\textbf{Power-type case (m=1).} Using the standard limit definition of the exponential function, we obtain:
\begin{equation}
\lim_{a \to 0} y = \lim_{a \to 0} (1 + axy^b)^{1/a} = e^{xy^b}.
\label{eq:exp_limit}
\end{equation}
Thus, the lucky equation at $a \to 0$ naturally transitions into the power-type generalization of the Lambert equation:
\begin{equation}
y = e^{xy^b} = M(1; 0; b; x), \qquad |x| < R.
\label{eq:lambert_power}
\end{equation}
The master numbers in this limit take the form:
\begin{equation}
N(1; 0; b; \ell) = \prod_{\gamma=1}^{\ell-1}(1 + b\ell) = (1 + b\ell)^{\ell-1}.
\label{eq:master_a0_power}
\end{equation}

\textbf{Logarithmic case (m=0).} Using the asymptotics $\ln(1+u) \approx u$ for small $u$, we obtain:
\begin{equation}
\lim_{a \to 0} w = \lim_{a \to 0} \frac{\ln(1 + ax e^{bw})}{a} = x e^{bw}.
\label{eq:lambert_log}
\end{equation}
Thus, the logarithmic form of the lucky equation at $a \to 0$ transitions into the classical Lambert equation:
\begin{equation}
w = x e^{bw} = M(0; 0; b; x), \qquad |x| < R.
\label{eq:lambert_log_eq}
\end{equation}
The master numbers in this limit take the form:
\begin{equation}
N(0; 0; b; \ell) = \prod_{\gamma=1}^{\ell-1}(b\ell) = (b\ell)^{\ell-1}.
\label{eq:master_a0_log}
\end{equation}

\textbf{Connection with the Lambert W function.} Of particular interest is the special case $b = -1$. Equation \eqref{eq:lambert_log_eq} takes the form:
\begin{equation}
w = x e^{-w}.
\label{eq:lambert_w}
\end{equation}
Multiplying both sides by $e^w$, we obtain the canonical defining equation for the Lambert $W$ function:
\begin{equation}
w e^w = x \quad \Longrightarrow \quad w = W(x).
\label{eq:lambert_w_canonical}
\end{equation}
Thus, we have strictly established the identity:
\begin{equation}
W_{0}(x) = M(0; 0; -1; x), \qquad |x| < R.
\label{eq:lambert_identity}
\end{equation}
This demonstrates that the Lambert $W$ function\cite{Corless1996, Valluri2000, Dubinov2006, Mezo2022} is not an isolated special function, but naturally embeds into the parametric family of logarithmic master series, as a special case of the lucky equation at $a \to 0$, $b = -1$.

Equally naturally, the Euler $T$ function\cite{Belkic2018, Belkic2019}, defined as the solution of $T e^{-T} = x$, also embeds into this family. In terms of the master series it takes the form
\begin{equation}
T(x) = M(0; 0; 1; x),
\label{eq:euler_T_master}
\end{equation}

\textbf{Special case a=0, b=0: classical functions.} With simultaneous zeroing of both parameters ($a \to 0$, $b = 0$), the lucky equation and its logarithmic form transition to the simplest classical functions:

\emph{Power-type case (m=1).} From formula \eqref{eq:exp_limit} at $b=0$, we obtain:
\begin{equation}
y = e^x = M(1; 0; 0; x).
\label{eq:exp_classical}
\end{equation}
The master numbers become identically equal to unity: $N(1; 0; 0; \ell) = 1^{\ell-1} = 1$, and the master series transforms into the classical exponential Taylor series:
\begin{equation}
e^x = 1 + x + \frac{x^2}{2!} + \frac{x^3}{3!} + \cdots
\label{eq:exp_taylor}
\end{equation}

\emph{Logarithmic case (m=0).} From formula \eqref{eq:lambert_log} at $b=0$, we obtain:
\begin{equation}
w = x = M(0; 0; 0; x).
\label{eq:identity_classical}
\end{equation}
This is the identity function --- the simplest possible master series.

Thus, at $a = 0$, $b = 0$, the lucky equation \emph{without any additional transformations} reduces to the classical exponential (and its logarithmic counterpart --- the identity function).

\textbf{Table of canonical master equations.}

\begin{table}[ht]
\centering
\caption{Canonical master equations}
\label{tab:canonical}
\begin{tabular}{cccc}
\toprule
$w = x$ & $w = x e^{b w}$ & $w = \frac{\ln(1 + a x)}{a}$ & $w = \frac{\ln\left( 1 + a x e^{b w} \right)}{a}$ \\
\midrule
$y = e^{x}$ & $y = e^{x y^{b}}$ & $y = (1 + a x)^{\frac{1}{a}}$ & $y = \left( 1 + a x y^{b} \right)^{\frac{1}{a}}$ \\
\midrule
$a = 0, b = 0$ & $a = 0, b \ne 0$ & $a \ne 0, b = 0$ & $a \ne 0, b \ne 0$ \\
\bottomrule
\end{tabular}
\end{table}

The principal branch of the solution to all canonical equations from this table, within the radius of convergence $|x|<R$, is given by a single master series $w = M(0; a; b; x)$, $y = M(1; a; b; x)$. All four equations of the lower row of the table ($m=1$) are unified by the fact that at $a = 0$, $b = 0$ their common master series reduces to the classical exponential $e^x$. In this sense, they represent a \textbf{two-parameter family of generalizations of the exponential}: the classical exponential is their common limiting point.

\subsection{Super-master series and the scaling identity}
\label{subsec:super_master}
The fundamental exponential identity \eqref{eq:exp_identity} $M(1;a;b;x) = \exp(M(0;a;b;x))$ gives rise to an important scaling property of master series, which we call the \textbf{super-master effect}.

Consider the master series at $m=0$:
\begin{equation}
M(0;a;b;x) = x + \frac{2b-a}{2!}x^2 + \frac{(3b-a)(3b-2a)}{3!}x^3 + \cdots
\label{eq:M0_explicit}
\end{equation}
Multiply it by an arbitrary constant $c$:
\begin{equation}
c \cdot M(0;a;b;x) = cx + c\frac{2b-a}{2!}x^2 + c\frac{(3b-a)(3b-2a)}{3!}x^3 + \cdots
\label{eq:cM0}
\end{equation}
The key observation: in the $\ell$-th term of the series, the number of factors in the master number is $\ell-1$, while the power of the argument is $\ell$. Therefore, under the substitution $a \to a/c$, $b \to b/c$, $x \to cx$, each factor of the master number acquires a common denominator $c$:
\begin{equation}
\left(\frac{b}{c}\ell - \frac{a}{c}\gamma\right) = \frac{1}{c}(b\ell - a\gamma),
\end{equation}
and the product of $\ell-1$ such factors yields a factor $c^{-(\ell-1)}$. Simultaneously, the argument $(cx)^\ell$ gives a factor $c^\ell$. As a result, each term of the series (except the constant term, which is absent at $m=0$) is multiplied by $c^\ell \cdot c^{-(\ell-1)} = c$, and we obtain the identity:
\begin{equation}
c \cdot M(0;a;b;x) = M\!\left(0;\frac{a}{c};\frac{b}{c};cx\right).
\label{eq:scaling_m0}
\end{equation}
Applying the fundamental exponential identity \eqref{eq:exp_identity} to both sides of \eqref{eq:scaling_m0}, we obtain the \textbf{power identity} for $m=1$:
\begin{equation}
M^c(1;a;b;x) = \exp\!\bigl(c \cdot M(0;a;b;x)\bigr) = \exp\!\left(M\!\left(0;\frac{a}{c};\frac{b}{c};cx\right)\right) = M\!\left(1;\frac{a}{c};\frac{b}{c};cx\right).
\label{eq:power_identity}
\end{equation}
Both identities --- \eqref{eq:scaling_m0} and \eqref{eq:power_identity} --- allow us to introduce a \textbf{fourth parameter} $c$ into the definition of the master series, giving rise to the \textbf{super-master series}~\cite{Berezin2026Py}:
\begin{equation}
S(m;a;b;x;c) = m + c\left(x + \frac{cm - a + 2b}{2!}x^2 + \frac{(cm - a + 3b)(cm - 2a + 3b)}{3!}x^3 + \cdots\right).
\label{eq:super_master}
\end{equation}

Thus, the super-master series $S$ is explicitly expressed in terms of the ordinary master series $M$ through the following identities (for $c \neq 0$):
\begin{align}
S(0; a; b; x; c) &= c \cdot M(0; a; b; x) = M\!\left(0; \frac{a}{c}; \frac{b}{c}; cx\right), \label{eq:S_m0_identity} \\
S(1; a; b; x; c) &= [M(1; a; b; x)]^c = M\!\left(1; \frac{a}{c}; \frac{b}{c}; cx\right). \label{eq:S_m1_identity}
\end{align}
These relations show that introducing the fourth parameter $c$ is equivalent either to scaling the logarithmic series or to raising the power series to a power, which in both cases reduces to a simple substitution of the argument and exponent parameters.

At $c=1$, the super-master series coincides with the ordinary one: $S(m;a;b;x;1) = M(m;a;b;x)$. At $m=0$, the parameter $c$ acts as a scaling factor; at $m=1$, it acts as a power exponent.

\textbf{Remark on the connection with the result of D\v{z}.~Belki\'c.} Formula (11.15) in \cite{Belkic2019} for an arbitrary power of a trinomial root of the equation $x=1+yx^\alpha$, obtained via Bell polynomials and the confluent Fox--Wright function, corresponds in the master-series terminology to the super-master $x^\beta=S(1; 1; \alpha; y; \beta)$ for the fixed parameter value $a=1$. Our super-master $x^c=S(1; a; b; y; c)$ generalizes this result to the case of arbitrary complex exponents $a, b \in \mathbb{C}$ in the equation $x^a=1+ayx^b$. This is necessary for the unambiguous identification of branches both inside and beyond the radius of convergence.

\paragraph{Connection with the series of Lambert and Euler.}
As shown in the review by Corless et al.~\cite{Corless1996}, Euler in 1779 obtained a series for the solution of the equation $y^A - y^B = (A-B)xy^{A+B}$ in the form
\[
y^c = 1 + cv + \frac{c(c+A+B)}{2!}v^2 + \frac{(c+A+2B)(c+2A+B)}{3!}v^3 + \cdots
\]
(formula 1.2 in~\cite{Corless1996}). Dividing Euler's equation by $y^B$ and making the substitution $a = A-B$, $b = A$ (so that $A = b$, $B = b-a$), we obtain the canonical master equation $y^a = 1 + axy^b$, and Euler's series for $y^c$ transforms into the super-master series $S(1; a; b; x; c)$.

Thus, in historical retrospect:
\begin{itemize}
\item the \textbf{master series} $M(0; a; b; x)$ generalizes the \emph{Belkić series} (2018) for the logarithm of the root of a trinomial equation;
\item the \textbf{master series} $M(1; a; b; x)$ generalizes the \emph{Lambert series} (1758) for the root of a trinomial equation;
\item the \textbf{super-master series} $S(1; a; b; x; c)$ coincides with the \emph{Euler series} (1779) for an arbitrary power of the root of a trinomial equation.
\end{itemize}

In fact, any of these series can be used as a unified tool for representing a wide class of functions: Euler's series is more general, and at $c=1$ the super-master series coincides with the master series. After establishing this universality, several natural questions arise: how to solve equations with a large number of terms using these series, how to find all roots of multivalued functions, and how to use the master series beyond its radius of convergence. To automate further research, the master series and super-master series have been implemented as functions in Python (see the repository~\cite{Berezin2026Py}).

The complete algorithm for computing the super-master series is given in Algorithm~\ref{alg:supermaster}. Readers who are not familiar with Python or Maple can implement it in any programming language by following the pseudocode.

\begin{algorithm}
\caption{SuperMaster: computation of the super-master series}
\label{alg:supermaster}
\begin{algorithmic}[1]
\Require $m, a, b, x \in \mathbb{C}$, $c \in \mathbb{C}$, $M \in \mathbb{N}$, $\varepsilon = 0.5 \times 10^{-D}$
\Ensure $(S, L)$ --- the value of the super-master series and the number of terms used
\State $f \gets c \cdot x$
\State $S \gets m + f$
\State $T \gets f$
\For{$L = 2$ \textbf{to} $M$}
    \State $N \gets x$
    \For{$g = 1$ \textbf{to} $L-1$}
        \State $N \gets N \cdot (c m - a g + b L) \cdot x / (g+1)$
    \EndFor
    \State $w \gets c \cdot N$
    \If{$w = 0$ \textbf{and} $T = 0$}
        \State \textbf{break}
    \EndIf
    \If{$w \neq 0$ \textbf{and} $|w/S| < \varepsilon$}
        \State \textbf{break}
    \EndIf
    \State $T \gets w$
    \State $S \gets S + w$
\EndFor
\State \Return $(S, L-1)$
\end{algorithmic}
\end{algorithm}

Here $D$ is the number of significant digits determined by the precision of the programming language used (for example, $D = 15$ for Python with double precision, $D = 10$ for Maple by default). The parameter $c$ can be any complex number; at $c = 1$, the super-master series coincides with the master series.

Algorithm~\ref{alg:supermaster} is written in a form convenient for implementation in languages without arbitrary precision (such as Python). In the Maple version, the computation of $f$ and $N$ is separated, which is possible due to symbolic computation. Mathematically, both variants are equivalent.

The returned value $L-1$ is the number of terms used in the series. It is necessary for studying the convergence rate and for subsequent accuracy checks in the \texttt{UltraRadical} function.

\subsection{General form of trinomial equations}
\label{subsec:general_trinomial}
\textbf{Algebraic trinomial.} Consider the general equation:
\begin{equation}
p\,Y^{a} = q + z\,Y^{b}, \qquad a, b \in \mathbb{C}, \quad p, q, z \in \mathbb{C}.
\label{eq:general_algebraic}
\end{equation}
Using the substitution
\begin{equation}
Y = y\left(\frac{q}{p}\right)^{1/a}, \qquad z = aqx\left(\frac{p}{q}\right)^{b/a},
\label{eq:substitution_algebraic}
\end{equation}
equation \eqref{eq:general_algebraic} reduces to the \textbf{canonical power form}:
\begin{equation}
y^{a} = 1 + axy^{b}.
\label{eq:canonical_from_general}
\end{equation}
The solution of the original equation is expressed through the master series:
\begin{equation}
Y_n = e^{f_n} \cdot M(1; a; b; x_n), \qquad x_n = \frac{z e^{bf_n}}{aq}, \qquad |x_n| < R.
\label{eq:solution_general_algebraic}
\end{equation}
\begin{equation}
f_n = \frac{\ln|q/p| + i\bigl[\arg(q/p) + 2\pi n\bigr]}{a}, \quad n \in \mathbb{Z}, \label{eq:f_def}
\end{equation}

Thus, we have established that the solutions of the multivalued algebraic trinomial equation of general form can be obtained via a \emph{single-valued} power master series $M(1; a; b; x)$, and the parameter $n \in \mathbb{Z}$ in formula \eqref{eq:f_def} ensures unambiguous identification of each solution branch.
\paragraph{Practical implementation.}
To compute the $n$-th branch of the solution of the equation $pY^a = q + zY^b$, the super-master series function is used. In a software implementation, this is written as follows:
\begin{equation}
Y_n = e^{f_n} \cdot \texttt{SuperMaster}(1, a, b, x_n, 1)[0],
\label{eq:practical_implementation}
\end{equation}
where $f_n$ is defined by formula \eqref{eq:f_def}, $x_n = \frac{z e^{b f_n}}{a q}$, and $[0]$ denotes the extraction of the series value from the pair $(S, L)$ returned by the function, where $S$ is the value of the super-master series and $L$ is the number of terms used.

\textbf{Remark on solution completeness.} The parametrization given above yields all solutions corresponding to \emph{one} of the three possible canonical transformations of the original trinomial (in this case --- the transformation where the dominant term is $p\,Y^a$). For equations with \textbf{complex} exponents, a complete description of all roots requires consideration of all three transformations; these will be detailed in Section~\ref{sec:three_transformations}. However, for \textbf{real} exponents all roots (all branches) are obtained through \emph{one} transformation --- under the condition $|x|<R$. Let us consider a characteristic example.

\textbf{Example: equation with fractional exponents.} Consider the equation
\begin{equation}
Y^{2/3} + 0.01\, Y^{1/2} + 1 = 0.
\label{eq:example_fractional}
\end{equation}
Reduction to canonical form yields parameters $a = 2/3$, $b = 1/2$, $p = 1$, $q = -1$, $z = -0.01$. Computing roots via formula \eqref{eq:solution_general_algebraic} for various $n \in \mathbb{Z}$ gives 4 formal solutions (since $b/a = 3/4$). However, equations with fractional, irrational, or complex exponents possess an important feature: many of their formal roots do not satisfy the original equation if exponentiation to powers $a$ and $b$ is performed via the principal branch of the complex logarithm.

Under standard verification via the principal branch, none of the obtained roots reduces equation \eqref{eq:example_fractional} to an identity. However, each root is strictly valid if the correct branch $u$ is used during substitution. The verification condition is written as
\begin{equation}
\exp\!\left(\frac{2}{3}\bigl[\ln|Y| + i(\arg Y + 2\pi u)\bigr]\right) + 0.01\,\exp\!\left(\frac{1}{2}\bigl[\ln|Y| + i(\arg Y + 2\pi u)\bigr]\right) + 1 = 0,
\label{eq:verification_example}
\end{equation}
where the parameter $u$ is applied identically to both power-law terms, despite the fact that only the exponent $a$ appears in definition~\eqref{eq:f_def}. The parameter $u$ is computed by the formula
\begin{equation}
u = \left\lceil \frac{\operatorname{Im}(f)}{2\pi} - \frac{1}{2} \right\rceil,
\label{eq:u_formula_example}
\end{equation}
where $f$ is defined by formula \eqref{eq:f_def}.

For this equation, all roots are obtained through a single transformation (with parameters $a = 2/3$, $b = 1/2$), since the convergence radii for the other two transformations do not admit convergent series for the given coefficients. Thus, in this case, the parametrization \eqref{eq:solution_general_algebraic} indeed gives the \emph{complete} solution of the equation.

\paragraph{The Root Passport of a Multinomial Equation.}
Thus, the roots of multinomial equations are not merely a disordered set of values. Each root possesses its own \textbf{certificate}, or \textbf{technical passport}, which, in addition to the value itself, specifies:
\begin{itemize}
\item $n$ --- the \textbf{branch index} of the multivalued function on which the root is evaluated;
\item $u$ --- the \textbf{verification index}. If the root appears in the equation raised to a fractional power, the result of this exponentiation is multivalued, and it is necessary to specify exactly which branch of the complex logarithm yields a valid identity.
\end{itemize}

In plasma physics, the arbitrary selection of ``any positive root'' is inadmissible; rather, the specific root that continuously extends from the equilibrium state is required. Without the parameter $n$, it is difficult to specify exactly which root is physically meaningful in the process under consideration. Without the parameter $u$, finding a root via iterative methods devolves into an exhaustive search across the sheets of the Riemann surface.

If all exponents are rational, each root possesses an entire family of equivalent branch and verification indices (for instance, for the equation $y^7 = 1 + Zy^3$, the indices $n = \dots, -7, 0, 7, 14, \dots$ yield the exact same root). However, this does not preclude unambiguously identifying the required root via its index $n$ and verifying it (or locating it iteratively) via its index $u$.

\textbf{Exponential trinomial.} Similarly, the exponential-type equation:
\begin{equation}
p\,Y = q\,e^{zY^{b}}
\label{eq:general_exponential}
\end{equation}
by the substitution $Y = y\frac{q}{p}$, $z = x\left(\frac{p}{q}\right)^{b}$ reduces to the canonical form $y = e^{xy^{b}}$.

Thus, the root of the general equation is expressed through the master series:
\begin{equation}
Y = \frac{q}{p} \cdot M\!\left(1; 0; b; z \left(\frac{q}{p}\right)^{b}\right).
\label{eq:general_exponential_solution}
\end{equation}

In the special case $b = 1$, this yields the Euler $T$ function. Indeed, the equation $T = x e^{T}$ corresponds to $a = 0$, $b = 1$, $m = 0$, so that
\begin{equation}
T(x) = M(0; 0; 1; x) = x \cdot M(1; 0; 1; x),
\label{eq:euler_T_via_master}
\end{equation}

\subsection{Radius of convergence of the master series}
\label{sec:convergence}

The master series $M(m; a; b; x) = m + \sum_{\ell=1}^{\infty} k_\ell x^\ell$, whose coefficients are given by
\begin{equation}
k_\ell = \frac{N_\ell}{\ell!} = \frac{1}{\ell!} \prod_{\gamma=1}^{\ell-1} (m - a\gamma + b\ell),
\label{eq:k_ell_def}
\end{equation}
has a finite radius of convergence $R$. According to the Cauchy--Hadamard formula:
\begin{equation}
\frac{1}{R} = \limsup_{\ell \to \infty} |k_\ell|^{\frac{1}{\ell}}.
\label{eq:cauchy_hadamard}
\end{equation}

To find the asymptotic behaviour as $\ell \to \infty$, we take the logarithm of $|k_\ell|$. The dominant contribution comes from the product $\prod_{\gamma=1}^{\ell-1}(b\ell - a\gamma)$, since $m$ is an additive constant. Factoring out $\ell$ from the product and replacing the sum of logarithms by an integral, using Stirling's formula $\ln(\ell!) \sim \ell\ln \ell - \ell$, we obtain:
\begin{equation}
\ln|k_\ell| \sim \ell \left( \int_{0}^{1} \ln|b - a \tau|\,d\tau + 1 \right) + o(\ell).
\label{eq:asymptotic_log}
\end{equation}

The integral can be evaluated analytically via the principal branch of the complex logarithm:
\begin{equation}
\int_{0}^{1} \ln|b - a \tau|\,d\tau = \operatorname{Re}\left( \frac{b\ln b - (b - a)\ln(b - a)}{a} \right) - 1.
\label{eq:integral_J}
\end{equation}

Consequently, the radius of convergence is:
\begin{equation}
R = \left| b^{-\frac{b}{a}} \cdot (b - a)^{\frac{b - a}{a}} \right|, \qquad a \neq 0,\; b \neq 0,\; b \neq a.
\label{eq:radius_final}
\end{equation}

For real parameters $a, b \in \mathbb{R}$ this formula simplifies to a compact form:
\begin{equation}
R = \frac{|1 - a/b|^{b/a}}{|b - a|}.
\label{eq:radius_compact}
\end{equation}

For degenerate cases the formula simplifies, and the case $b=a$ is elegantly resolved via the master identities:
\begin{itemize}
    \item \textbf{When $b = 0$:} $R = \frac{1}{|a|}$ (the classical binomial series).
    \item \textbf{When $a \to 0$:} $R = \frac{1}{|b| e}$ (the generalised exponential / Lambert function).
    \item \textbf{When $b = a$:} The equation takes the form $y^a = 1 + axy^a$, which is algebraically equivalent to $y = (1 - ax)^{-1/a}$. By the master identity $M(m;a;b;x) = M(m;-a;b-a;x)$, this corresponds precisely to the case $b=0$ with the parameter substitution $a \to -a$. Therefore, the radius of convergence is $R = \frac{1}{|a|}$.
\end{itemize}

\subsubsection*{Behaviour on the boundary of convergence}
The behaviour of the power series on the boundary $|x| = R$ depends on the specific equation. For example:

\begin{itemize}
\item The series for $\sqrt{1+x}$ converges for $|x| \leq 1$.
\item The series for $\frac{1}{\sqrt{1+x}}$ converges for $-1 < x \leq 1$.
\end{itemize}

For brevity, in this work we use the unified convergence condition $|x| < R$. It should be understood that the actual behaviour depends on the specific equation:

\begin{itemize}
\item For some equations: $x < R$.
\item For other equations: $x \leq R$.
\end{itemize}

The derivatives and integrals of master series share the same radius of convergence $R$ as the original series; the only possible difference is whether the series converges strictly inside the disk $|x| < R$ or admits convergence on the boundary $|x| = R$ itself.

\subsection{Analytic Continuations}
\label{app:analytic_continuation}

The convergence radius of the exponential, sine and cosine is infinite. Binomial series has limited number of terms, since all subsequent master series terms equal zero. Therefore convergence question for this series is irrelevant.
\[
(1 + x)^{a} = M(1;1/a;0;ax)\quad ,\quad a \in \mathbb{N}
\]
\[
M(1;1/2;0;2x) = 1 + 2x + \frac{4x^{2}}{2}\left( 1 - \frac{1}{2} \right) + \frac{8x^{3}}{3!}\left( 1 - \frac{1}{2} \right)\left( 1 - \frac{2}{2} \right) = 1 + 2x + x^{2}
\]
\[
M(1;1/3;0;3x) = 1 + 3x + \frac{9x^{2}}{2}\left( 1 - \frac{1}{3} \right) + \frac{27x^{3}}{3!}\left( 1 - \frac{1}{3} \right)\left( 1 - \frac{2}{3} \right) + 0 + 0 + \cdots
\]

For multivalued functions, the question of limited convergence radius is especially important. Each branch of a multivalued function is unique, and to identify it we refer to the power master series with one set of parameters $M(m; a; b; x)$ --- this is the so-called \emph{origin} of the branch, its local representation in a neighbourhood of some point. However, for multivalued functions this initial series has a finite radius of convergence $R$. The same branch is defined in another region (beyond $|x| < R$) by the same universal power master series, but now with different parameters $M(m; \tilde{a}; \tilde{b}; \tilde{x})$, convergent in that region. We call such a series the \textbf{analytic continuation} of the specifically indicated branch.

\subsubsection*{Analytic continuations for $b=0$}
\label{app:analytic_continuation_b0}

For $b=0$, the convergence condition \eqref{eq:kernel_convergence} ($|K|<1$) is not required --- it is sufficient to satisfy condition \eqref{eq:binomial_convergence} ($|ax|<1$), as shown in the derivation of the lucky series. This means that the kernel modulus $|K|$ can be arbitrarily large, and the question of analytic continuation for the generalized logarithm and binomial root is resolved straightforwardly.

Consider the function $x^{1/a}$ for arbitrary $|x|$. For $|x| < 1$, the standard representation with parameter $s=1$ is used. For $|x| \geq 1$, we apply the identity $|x|^{1/a} = |x^{-1}|^{-1/a}$, which is equivalent to replacing the parameter $a \to -a$ and the argument $x \to x^{-1}$. In both cases, the factor selecting the specific branch of the multivalued function remains unchanged:

\begin{equation}
x^{\frac{1}{a}} = M\left( 1; sa; 0; \frac{|x^{s}| - 1}{sa} \right) 
e^{\frac{(\arg(x) + 2\pi n)i}{a}}, \quad n \in \mathbb{Z},
\label{eq:analytic_cont_b0}
\end{equation}
where the parameter $s$ is determined by the convergence condition:
\begin{equation}
s = \begin{cases}
1, & |x| < 1, \\
-1, & |x| \geq 1.
\end{cases}
\label{eq:h_parameter}
\end{equation}

Thus, the analytic continuation of the function $x^{1/a}$ beyond the circle $|x| < 1$ is realized by switching the parameter $s$ and the corresponding change of the master series argument, while the series structure and branch factor remain preserved.

An analogous formula ensures the continuity of any branch $n$ of the generalized complex logarithm:
\begin{equation}
\frac{\ln(x)}{a} = M\left( 0; sa; 0; \frac{|x^{s}| - 1}{sa} \right) 
+ \frac{(\arg(x) + 2\pi n)i}{a}, \quad n \in \mathbb{Z}.
\label{eq:analytic_cont_log_b0}
\end{equation}
\paragraph*{Conclusion for the case $b=0$.}
Thus, for $b=0$ and $a \neq 0$, the multivaluedness of the solution is entirely captured by a single integer parameter $n \in \mathbb{Z}$, which indexes the branches of the complex root or the generalized logarithm. Analytic continuation in this case is trivial: it suffices to replace $a \to -a$ and $x \to x^{-1}$ while keeping the index $n$ unchanged. Each branch is uniquely identified by its index $n$ throughout the entire complex plane. Moreover, the master series $M(m; a; b; x)$ itself does not depend on the value of $n$: the entire multivaluedness of the solution is concentrated exclusively in the exponential factor $e^{(\arg(x) + 2\pi n)i/a}$, while the series remains a single universal function for all branches.

\begin{remark}[Computational nature of the discrete parameter $s$]
A natural question arises: why does the parameter $s$ in the analytic continuation formulas take only discrete values $s = \pm 1$, instead of being a smooth function, such as $s = 1 - |x|$ or $s = -\ln|x|$, which would ensure a continuous matching of the series?

From the perspective of formal algebra, the identity $(x^s)^{1/(sa)} = x^{1/a}$ holds for any $s \neq 0$. However, a computational complexity analysis renders continuous values of $s$ meaningless. To compute the new argument of the master series $X_{\text{new}} = x^s$ for a continuous $s$, one must perform exponentiation to the power of $s$, which is computationally equivalent to the original root-extraction problem $x^{1/a}$. We thus fall into a computational tautology: to evaluate the series, we must first compute the series.

For the analytic continuation method to be practically meaningful, the argument transformation $x \to x^s$ must be computationally trivial (requiring $O(1)$ operations). There exist exactly two such transformations that change the modulus of a number without resorting to power series: the identity transformation ($s = 1$, $x \to x$) and inversion ($s = -1$, $x \to 1/x$).

Thus, the discrete choice $s = \operatorname{sgn}(1 - |x|)$ is dictated not by any incompleteness of the method, but by \textbf{computational complexity constraints}.
\end{remark}

\subsubsection*{Branch indices $n, h, k$ for $a, b \neq 0$, $a \neq b$, and $a, b \in \mathbb{R}$}
\label{sec:three_transformations}

The lucky equation has nonzero and unequal parameters $a$ and $b$. Under these parameters, the kernel modulus $|K|$ is bounded by unity (see condition \eqref{eq:kernel_convergence}), which ensures the convergence of the binomial expansion $(1+K)^{\ell b}$. Analytic continuation of specific branches in this case is more complex than for $b=0$.

Consider the original trinomial equation
\begin{equation}
S_2 Y^{\beta_2} + S_1 Y^{\beta_1} + S_0 = 0,
\label{eq:trinomial_general}
\end{equation}
where the coefficients $S_0, S_1, S_2 \in \mathbb{C}$ are complex and the exponents $\beta_2, \beta_1 \in \mathbb{R}$ are real, with the ordering
\begin{equation}
\beta_2 > \beta_1 > 0
\label{eq:beta_ordering}
\end{equation} adopted to ensure consistent branch indexing.

There are 6 ways to transform equation \eqref{eq:trinomial_general} to the canonical form $p_{j}Y^{a_{j}} = q_{j} + z_{j}Y^{b_{j}}$. The transformation names reflect the fraction $q_{j}/p_{j}$, from which the parameter $f_{j,N_{j}}$ is formed:
\begin{equation}
f_{j,N_{j}} = \frac{\ln|q_{j}/p_{j}| + i\bigl(\arg(q_{j}/p_{j}) + 2\pi N_{j}\bigr)}{a_{j}}.
\label{eq:f_definition}
\end{equation}

In practice, three transformations suffice $S_0/S_2$, $S_2/S_1$, $S_1/S_0$ ($j=1,2,3$):

\begin{table}[ht]
\centering
\caption{Three canonical transformations of the trinomial equation}
\label{tab:three_transformations}
\begin{tabular}{ccccccccc}
\toprule
j & Transformation & Form & $p_{j}$ & $q_{j}$ & $a_{j}$ & $b_{j}$ & $z_{j}$ & $N_{j}$ \\
\midrule
1 & $S_0/S_2$ & $S_2 Y^{\beta_2} = -S_0 - S_1 Y^{\beta_1}$ & $S_2$ & $-S_0$ & $\beta_2$ & $\beta_1$ & $-S_1$ & $n$ \\
\midrule
2 & $S_2/S_1$ & $S_1 Y^{\beta_1-\beta_2} = -S_2 - S_0 Y^{-\beta_2}$ & $S_1$ & $-S_2$ & $\beta_1-\beta_2$ & $-\beta_2$ & $-S_0$ & $h$ \\
\midrule
3 & $S_1/S_0$ & $S_0 Y^{-\beta_1} = -S_1 - S_2 Y^{\beta_2-\beta_1}$ & $S_0$ & $-S_1$ & $-\beta_1$ & $\beta_2-\beta_1$ & $-S_2$ & $k$ \\
\bottomrule
\end{tabular}
\end{table}

The mirror reflections of these transformations ($S_2/S_0$, $S_1/S_2$, $S_0/S_1$) yield the same roots and are not used in this work. Here we present the table of mirror transformations to illustrate the identity $M(m;a;b;x) = M(m;-a;b-a;x)$, which will be discussed in the study of master numbers.

\begin{table}[ht]
\centering
\caption{Mirror transformations of the trinomial equation}
\label{tab:three_transformations_2}
\begin{tabular}{ccccccccc}
\toprule
j & Transformation & Form & $p_{j}$ & $q_{j}$ & $a_{j}$ & $b_{j}$ & $z_{j}$ & $N_{j}$ \\
\midrule
4 & $S_2/S_0$ & $S_0 Y^{-\beta_2} = -S_2 - S_1 Y^{\beta_1-\beta_2}$ & $S_0$ & $-S_2$ & $-\beta_2$ & $\beta_1-\beta_2$ & $-S_1$ & $-n$ \\
\midrule
5 & $S_1/S_2$ & $S_2 Y^{\beta_2-\beta_1} = -S_1 - S_0 Y^{-\beta_1}$ & $S_2$ & $-S_1$ & $\beta_2-\beta_1$ & $-\beta_1$ & $-S_0$ & $-h$ \\
\midrule
6 & $S_0/S_1$ & $S_1 Y^{\beta_1} = -S_0 - S_2 Y^{\beta_2}$ & $S_1$ & $-S_0$ & $\beta_1$ & $\beta_2$ & $-S_2$ & $-k$ \\
\bottomrule
\end{tabular}
\end{table}
To obtain the roots of the original equation \eqref{eq:trinomial_general}, we use the parameters $p_{j}$, $q_{j}$, $a_{j}$, $b_{j}$, $z_{j}$, $N_{j}$ from table \ref{tab:three_transformations} to determine the master series argument:
\begin{equation}
x_{j,N_j} = \frac{z_{j} e^{b_{j}f_{j,N_j}}}{a_{j}q_{j}}, \quad \text{where } f_{j,N_j} \text{ is given by formula } \eqref{eq:f_definition}.
\label{eq:x_argument}
\end{equation}
Then we find the root of branch with index $j,N_j$ by the formula:
\begin{equation}
Y_{j,N_j} = e^{f_{j,N_j}} \cdot M(1;a_{j};b_{j};x_{j,N_j}).
\label{eq:root_formula}
\end{equation}

However, the use of a two-dimensional index $(j, N_j)$ defeats the very purpose of branch identification: the transformation number $j$ is merely a technical parameter, whereas the physical meaning resides in a single identifier of a continuous branch across the entire Riemann surface. Therefore, the key task is to establish the correspondence: which branch $h$ (transformation $j=2$) or $k$ (transformation $j=3$) is the analytic continuation of the given branch $n$ (transformation $j=1$). Once this correspondence is established, all three transformations describe a single continuous branch, uniquely identified by the index $n$.

The complex ultra-radical ($a,b,x \in \mathbb{C}$) differs qualitatively from the ultra-radical operating with real exponents. For complex exponents, convergence at $|x| < R$ depends on the value of $N$.

The real ultra-radical ($a,b \in \mathbb{R}$, $x \in \mathbb{C}$) is simpler: either all branches are determined through the $S_0/S_2$ transformation and converge for any $n$, or all roots are found through the $S_2/S_1$ and $S_1/S_0$ transformations and converge for any integer values of $h$ and $k$, respectively.

This property has a rigorous algebraic justification. For real exponents, the moduli of the arguments $|x_j|$ are independent of the chosen branch indices $n$, $h$, $k$. Moreover, the conditions for reaching the convergence boundary are identical for all three transformations:
\[
|x_1| = R_1 \iff |x_2| = R_2 \iff |x_3| = R_3,
\]
However, since $x_1$ depends linearly on the parameter $S_1$, while $x_2$ and $x_3$ are power-law functions with negative exponents (for $\beta_1 < \beta_2$), the inequalities point in opposite directions: on one side of the critical point, the base series ($j=1$) converges, while on the other side, the conjugate series ($j=2$ and $j=3$) converge.

Here arises the key question: which of the branches $h$ or $k$ is the analytic continuation of the given branch $n$, if the series for the $S_0/S_2$ transformation diverges? The answer to this question is provided by the geometric criterion, discussed below.

\subsubsection*{Geometric criterion for selecting the conjugate series of the real ultra-radical}

The real ultra-radical $U_n(\beta_2, \beta_1; Z/\beta_2)$, which is studied in detail in~\cite{Berezin2025_arXiv}, is defined as the $n$-th branch of the solution to the equation
\begin{equation}
Y^{\beta_2} = 1 + Z Y^{\beta_1},
\label{eq:real_ultraradical_def}
\end{equation}
where the parameter $Z \in \mathbb{C}$, and the exponents $\beta_2, \beta_1 \in \mathbb{R}$ satisfy the condition $\beta_2 > \beta_1 > 0$.

To apply this definition to the general equation~\eqref{eq:trinomial_general}, it must first be reduced to the unit form~\eqref{eq:real_ultraradical_def} via the substitution~\eqref{eq:substitution_algebraic}, which normalizes the coefficient of the leading term to unity and the constant term to $-1$.

Consider the equation $y^7 = 1 + Z y^3$. We will vary the parameter $Z$ from $-\infty$ to $+\infty$ and monitor the continuity of the graph of the branch with index $n=0$.

When $|Z/7| < R$, the base $S_0/S_2$ transformation is used. The base points $e^f$ form a regular heptagon on the complex plane (the 7th roots of unity). Since the exponent 3 is an integer, exactly one root of the original equation corresponds to each base point. All 7 roots are obtained through a single transformation; the branch index $n=0$ is specified by the user, and the algorithm substitutes $n, 7, 3, Z$ into formula~\eqref{eq:root_formula} according to Table~\ref{tab:three_transformations}, yielding the root of branch 0. Since 7 and 3 are integers, the same root will be given by indices $n = \dots, -7, 0, 7, 14, \dots$. The smaller $|Z|$ is, the closer the roots are to their base points.

However, when $|Z/7| > R$, two different transformations are used: $S_2/S_1$ and $S_1/S_0$. In the $S_2/S_1$ transformation, the base points form a square $(1/Z)^{1/(3-7)}$ on the complex plane, while in $S_1/S_0$ they form an equilateral triangle $(-Z/1)^{-1/3}$. \textbf{The larger $|Z|$ is, the closer the roots of the original equation are to their base points --- the vertices of the square and the triangle.} The total number of roots remains unchanged (still 7), but the periodicity of identical roots is specific to each transformation. If one of the exponents is irrational, there is no periodicity of identical roots: each value of the indices $n, h, k$ yields a unique root. Searching beyond the radius of convergence for a branch continuous with index $n=0$ by trial and error is computationally infeasible.

The real ultra-radical algorithm~\cite{Berezin2026Py} solves this problem for an arbitrary branch index $n \in \mathbb{Z}$ as follows. The concept of \emph{candidates} is introduced --- the centers of sectors on the imaginary axis:
\begin{equation}
L_{1,n} = \operatorname{Im}(\beta_1 f_{1,n}), \quad
L_{2,h} = \operatorname{Im}(\beta_1 f_{2,h}), \quad
L_{3,k} = \operatorname{Im}(\beta_1 f_{3,k}),
\label{eq:sector_centers}
\end{equation}
where $\beta_1$ is the exponent from the original equation (in our example, $\beta_1=3$). \textbf{Important note:} the sector is multiplied precisely by $\beta_1$ of the original equation, not by its analogue under the transformation. That is, the true length of any sector is multiplied by a common factor. This is necessary in cases where the exponents are fractional (for example, in the equation $y^{2/3} + 100 y^{1/2} + 1 = 0$).

For real $\beta_2, \beta_1$, the expressions~\eqref{eq:sector_centers} simplify to:
\begin{equation}
L_{1,n} = \frac{2\pi \beta_1 n}{\beta_2}, \quad
L_{2,h} = \frac{\beta_1[\arg(1/Z) + 2\pi h]}{\beta_1 - \beta_2}, \quad
L_{3,k} = \frac{\beta_1[\arg(-Z) + 2\pi k]}{-\beta_1}.
\label{eq:sector_centers_real}
\end{equation}

The centers $L_{1,n}$ partition the imaginary axis into sectors $[g_n, G_n] = [L_{1,n} - \pi \beta_1/\beta_2,\ L_{1,n} + \pi \beta_1/\beta_2]$. The frequency of points $L_{1,n}$ on the imaginary axis equals the sum of the frequencies of points $L_{2,h}$ and $L_{3,k}$:
\[
\frac{\beta_2}{2\pi\beta_1} = \frac{\beta_2-\beta_1}{2\pi\beta_1} + \frac{\beta_1}{2\pi\beta_1}.
\]
This relation reflects the fact that, under the condition $\beta_2 > \beta_1 > 0$, the points $L_{2,h}$ and $L_{3,k}$ uniformly ``fill'' the sectors formed by the points $L_{1,n}$, without overlapping each other within a single sector. Therefore, each sector $n$ can contain at most one integer value of $h$ and at most one integer value of $k$. Moreover, due to the mutual arrangement of the base points $h = 0$ and $k = 0$, these integer values cannot simultaneously belong to the same sector $n$ (except for the degenerate cases at the boundaries considered below). If one of the points $L_{2,h}$ or $L_{3,k}$ lies inside sector $n$, then the corresponding value of $h$ or $k$ is the analytic continuation of the specified branch $n$.

Algorithmically, this is implemented as follows. The boundaries of sector $n$ are mapped into the coordinate system of the conjugate transformation, yielding real values $h_{\min}$ and $h_{\max}$ (which, in general, may be fractional). If the interval $[h_{\min}, h_{\max}]$ contains an integer, this integer is adopted as the required index of the branch $h$ conjugate to branch $n$. Note that under the condition $\beta_2 > \beta_1 > 0$, there cannot be two integer values of $h$ in this range. The range of fractional values of $h$ is found from the equations $L_{2,h_{\min}} = L_{1,n} - \pi \beta_1/\beta_2$ and $L_{2,h_{\max}} = L_{1,n} + \pi \beta_1/\beta_2$, solved for $h$. Similarly, the presence of an integer $k$ in the same sector $n$ is determined.

In such cases, the choice of candidate is unambiguous. However, there are cases where within the specified sector $n$ there are no candidates with integer values of $h$ or $k$, but on the boundary between sectors $(n-1, n)$ or $(n, n+1)$ there are two candidates $h$ and $k$ simultaneously. From the perspective of branch continuity, there is no difference in choosing either as the analytic continuation, since at $|Z/\beta_2|=R$ both branches intersect.

In this degenerate case at sector boundaries, where multiple candidates for continuation exist, we propose the following convention: when two branches $h$ and $k$ intersect on the boundary of sectors $n-1$ and $n$, the branch with index $n$ is continued via the $h$-series, and the branch with index $n-1$ via the $k$-series. This ensures synchronization of the numbering of branches $0$ and $-1$ of the ultra-radical with branches $W_0$ and $W_{-1}$ of the Lambert $W$-function. A rigorous proof of the sector-based conjugation criterion and its extension to the ultralogarithm and complex exponents are presented in~\cite{Berezin2026_ULog}.

When plotting $Y_0(Z)$ for equation~\eqref{eq:real_ultraradical_def} as $Z$ varies from $-\infty$ to $+\infty$, passing through the point $Z=0$ results in a phase shift of $\pi$ in the argument of $Z$. Consequently, the choice between $h=0$ and $k=0$ instead of $n=0$ when crossing the boundary $|Z/\beta_2|=R$ depends directly on the sign of $Z$. This criterion is universal for any $n$: the selection of the conjugate series $h$ or $k$ upon crossing $|Z/\beta_2|=R$ is determined not only by the branch index $n$, but also by the argument $\arg(Z)$.

Thus, the geometric criterion transforms the index $n$ into a true one-dimensional branch identifier. It preserves its continuous mathematical and physical meaning across the entire complex plane, uniquely labeling the solution regardless of which transformation ($n$, $h$, or $k$) is currently active.

\subsection{The Real Master-Join Ultra-Radical}
\label{sec:real_ultraradical}

Let $y = \sqrt[{n;a;b}]{x}$ denote the \textbf{real Master-Join ultra-radical} --- the $n$-th branch of the solution to the equation
\begin{equation}
y^{a} = 1 + a x y^{b},
\label{eq:ultraradical_def}
\end{equation}
where $a, b \in \mathbb{R}$ and $x \in \mathbb{C}$.

The complete algorithm for computing the real ultra-radical $U_n(a;b;x) \equiv \sqrt[{n;a;b}]{x}$, including the geometric criterion for selecting the conjugate series beyond the radius of convergence, is formalized as a deterministic computational procedure.

For convenience of description, we introduce an auxiliary function $\text{SM}(j, N)$, which computes the value of the master series for the $j$-th transformation ($j \in \{1, 2, 3\}$) and branch $N \in \mathbb{Z}$:
\begin{equation}
\begin{aligned}
f &= \frac{\ln|q_j/p_j| + i\bigl(\arg(q_j/p_j) + 2\pi N\bigr)}{a_j}, \\
u &= \left\lceil \frac{\operatorname{Im}(f)}{2\pi} - \frac{1}{2} \right\rceil, \quad v = e^f, \quad V = e^{b_j f}, \\
[U, \ell] &= \text{SuperMaster}\left(1, a_j, b_j, \frac{z_j \cdot V}{a_j q_j}\right),
\end{aligned}
\label{eq:sm_function}
\end{equation}
The function $\text{SM}(j, N)$ returns the tuple $[v \cdot U,\; \ell,\; [j, N],\; u,\; f]$, where $v \cdot U$ is the root, $\ell$ is the number of master-series terms used, $[j, N]$ identifies the transformation and branch, $u$ is the verification index, and $f$ is the auxiliary quantity preserved for analytical studies.

The complete procedure for computing the ultra-radical is presented in Algorithm~\ref{alg:ultraradical_real}.

\begin{algorithm}[H]
\caption{Initialization of transformations}
\label{alg:init_transforms}
\begin{algorithmic}[1]
\Require $a_1, b_1 \in \mathbb{R}$, $x \in \mathbb{C}$
\Ensure Arrays $p, q, z, a, b$ for three transformations
\State $S_2 \gets 1, \quad S_1 \gets -a_1 x, \quad S_0 \gets -1$
\State $p \gets [S_2, S_1, S_0], \quad q \gets [-S_0, -S_2, -S_1]$
\State $z \gets [-S_1, -S_0, -S_2]$
\State $a \gets [a_1, b_1-a_1, -b_1], \quad b\gets [b_1, -a_1, a_1-b_1]$
\State \Return $p, q, z, a, b$
\end{algorithmic}
\end{algorithm}

\begin{algorithm}[H]
\caption{Computation of sector boundaries}
\label{alg:sector_bounds}
\begin{algorithmic}[1]
\Require $n \in \mathbb{Z}$, arrays $a, b$ from Algorithm~\ref{alg:init_transforms}, $q, p$
\Ensure Boundaries $[h_{\min}, h_{\max}]$ and $[k_{\min}, k_{\max}]$
\State $h_{m1} \gets \frac{(n-1/2) \cdot a_2}{a_1} - \arg(q_2/p_2)/2\pi$
\State $h_{m2} \gets \frac{(n+1/2) \cdot a_2}{a_1} - \arg(q_2/p_2)/2\pi$
\State $k_{m1} \gets \frac{(n-1/2) \cdot a_3}{a_1} - \arg(q_3/p_3)/2\pi$
\State $k_{m2} \gets \frac{(n+1/2) \cdot a_3}{a_1} - \arg(q_3/p_3)/2\pi$
\State $h_{\min} \gets \min(h_{m1}, h_{m2}), \quad h_{\max} \gets \max(h_{m1}, h_{m2})$
\State $k_{\min} \gets \min(k_{m1}, k_{m2}), \quad k_{\max} \gets \max(k_{m1}, k_{m2})$
\State \Return $h_{\min}, h_{\max}, k_{\min}, k_{\max}$
\end{algorithmic}
\end{algorithm}

\begin{algorithm}[H]
\caption{Computation of the real ultra-radical $U_n(a; b; x)$}
\label{alg:ultraradical_real}
\begin{algorithmic}[1]
\Require $a_1, b_1 \in \mathbb{R}$ ($a_1 > b_1 > 0$), $x \in \mathbb{C}$, $n \in \mathbb{Z}$
\Ensure $[y, \ell, \text{branch\_info}, u, f]$
\State $[p, q, z, a, b] \gets \text{InitTransforms}(a_1, b_1, x)$ \Comment{Algorithm~\ref{alg:init_transforms}}
\State $[h_{\min}, h_{\max}, k_{\min}, k_{\max}] \gets \text{SectorBounds}(n, a, b, q, p)$ \Comment{Algorithm~\ref{alg:sector_bounds}}
\State $h \gets$ integer in the range $[h_{\min}, h_{\max}]$, otherwise $h \gets \text{!}$
\State $k \gets$ integer in the range $[k_{\min}, k_{\max}]$, otherwise $k \gets \text{!}$
\State $R \gets \left| \frac{(1 - a_1/b_1)^{b_1/a_1}}{a_1 - b_1} \right|$
\If{$|x| < R$}
    \State \Return result of the call $\text{SM}(1, n)$
\EndIf
\If{$k = \text{!}$}
    \State \Return result of the call $\text{SM}(2, h)$ \Comment{The only candidate is $h$}
\ElsIf{$h = \text{!}$}
    \State \Return result of the call $\text{SM}(3, k)$ \Comment{The only candidate is $k$}
\Else
    \If{$\bigl(\arg(q_2/p_2) + 2\pi h\bigr) \frac{b_1}{a_2} < 2\pi n \frac{b_1}{a_1}$}
        \State \Return result of the call $\text{SM}(2, h)$ \Comment{Intersection to the left of the center of sector $n$}
    \Else
        \State \Return result of the call $\text{SM}(3, k)$ \Comment{Intersection to the right of the center of sector $n$}
    \EndIf
\EndIf
\end{algorithmic}
\end{algorithm}

Thus, the analytical formula for the real ultra-radical can be written compactly:
\begin{equation}
\sqrt[{n;a;b}]{x} = v_{j,N} \cdot M(1; a_j; b_j; x_{j,N}),
\label{eq:ultraradical_compact}
\end{equation}
where the transformation index $j \in \{1, 2, 3\}$ and the branch index $N \in \{n, h, k\}$ are determined by Algorithm~\ref{alg:ultraradical_real}, and the parameters are given by:
\begin{equation}
v_{j,N} = e^{f_{j,N}}, \quad
f_{j,N} = \frac{\ln|q_j/p_j| + i\bigl(\arg(q_j/p_j) + 2\pi N\bigr)}{a_j}, \quad
x_{j,N} = \frac{z_j \, e^{b_j f_{j,N}}}{a_j q_j}.
\label{eq:ultraradical_params}
\end{equation}

\subsubsection*{Derivative of the Ultra-Radical}

Differentiating the identity~\eqref{eq:ultraradical_def} with respect to $x$, we obtain:
\[
a y^{a-1} y' = a y^{b} + a b x y^{b-1} y'.
\]
Dividing by $a$ (for $a \neq 0$) and grouping the terms containing $y'$, we find:
\[
y' \left( y^{a-1} - b x y^{b-1} \right) = y^b.
\]
Dividing the numerator and denominator by $y^{b-1}$, we arrive at a compact and computationally convenient form:
\begin{equation}
\boxed{\,
\frac{d}{dx}\sqrt[{n;a;b}]{x} 
= \frac{y^{b-a+1}}{1 - b x y^{b-a}}
= \frac{y}{y^{a-b} - b x}.
\,}
\label{eq:ultra_derivative}
\end{equation}
Both expressions in~\eqref{eq:ultra_derivative} are algebraically identical; the first emphasizes the degree of deviation from unity, while the second is convenient when substituting the original equation $y^a = 1 + axy^b$.

Consequently, the ultra-radical is a solution to a first-order nonlinear ordinary differential equation:
\begin{equation}
\label{eq:ultra_ode}
\boxed{\,
y' = \frac{y^{b-a+1}}{1 - b x y^{b-a}},
\qquad y(0)=1.
\,}
\end{equation}
Conversely, any ODE reducible to the form~\eqref{eq:ultra_ode} has a solution expressible through the ultra-radical.

\subsubsection*{Integrals of the Ultra-Radical}

To compute the integral $\int y\,dx$, we use the method of integration via the inverse function. From the defining equation $y^a = 1 + axy^b$, we express $x$ as a function of $y$:
\begin{equation}
x = \frac{y^{a-b} - y^{-b}}{a}.
\label{eq:x_from_y}
\end{equation}

Applying the formula $\int y\,dx = xy - \int x\,dy$ and substituting~\eqref{eq:x_from_y}, we find:
\begin{align*}
\int y\,dx &= xy - \frac{1}{a}\int \bigl( y^{a-b} - y^{-b} \bigr)\,dy \\
&= xy - \frac{1}{a}\left[ \frac{y^{a-b+1}}{a-b+1} - \frac{y^{1-b}}{1-b} \right] + C,
\end{align*}
where it is assumed that $a-b+1 \neq 0$ and $b \neq 1$.

Thus, the general formula for the indefinite integral takes the form:
\begin{equation}
\boxed{\,
\int \sqrt[{n;a;b}]{x}\,dx 
= x\sqrt[{n;a;b}]{x} 
- \frac{1}{a}\left[ 
\frac{(\sqrt[{n;a;b}]{x})^{a-b+1}}{a-b+1} 
- \frac{(\sqrt[{n;a;b}]{x})^{1-b}}{1-b} 
\right] + C.
\label{eq:ultra_integral_general}
\,}
\end{equation}

\textbf{Special case $b=1$.} When $b=1$, formula~\eqref{eq:ultra_integral_general} takes the form:
\begin{equation}
\int \sqrt[{n;a;1}]{x}\,dx 
= x\sqrt[{n;a;1}]{x} 
- \frac{1}{a}\left[ 
\frac{(\sqrt[{n;a;1}]{x})^{a}}{a} 
- \ln\!\bigl( \sqrt[{n;a;1}]{x} \bigr) 
\right] + C.
\label{eq:ultra_integral_b1}
\end{equation}

\subsubsection{Applications of the Real Ultra-Radical}
\label{sec:applications}
In many stationary physical processes, the key principle is the conservation law for the total energy per particle:
\begin{equation}
E_{\text{kin}} + E_{\text{int}} + E_{\text{pot}} = \text{const}.
\label{eq:energy_const}
\end{equation}
If we consider two states of the system (initial state $1$ and current state $2$) and subtract their energy balances, the constant disappears, and we obtain the sum of increments equal to zero:
\begin{equation}
\Delta E_{\text{kin}} + \Delta E_{\text{int}} + \Delta E_{\text{pot}} = 0.
\label{eq:energy_delta}
\end{equation}
For a wide class of media, kinetic and internal energies depend on characteristic state variables (velocity, density, temperature) in a power-law manner. After subtracting the initial state, each increment takes the form $A(1 - \alpha)$, where $\alpha$ is the ratio of the current value to the initial one. Introducing a universal dimensionless variable $w$ that relates the ratios of these quantities through the continuity equation or equation of state, equation~\eqref{eq:energy_delta} reduces to the canonical form:
\begin{equation}
\mathcal{K}\left(1 - \frac{1}{w^2}\right) + \frac{\Theta}{\gamma-1}\left(1 - w^{\gamma-1}\right) - \Phi = 0.
\label{eq:general_trinomial}
\end{equation}
Here $\mathcal{K}$, $\Theta$, $\Phi$ are energy scales depending on the physics of the process, $\gamma$ is the adiabatic index, and $w$ is a universal unknown describing the degree of deviation of the system from equilibrium. In the considered applications, $w$ represents the ratio of physical quantities (densities, velocities, concentrations), so by definition $w > 0$. Multiplying~\eqref{eq:general_trinomial} by $w^2$ immediately leads to a trinomial equation solvable via the ultraradical. For such equations, the ultraradical always returns positive real solutions~\cite{Berezin2026Py} on branches $n=-1$ and $n=0$. The physically realizable branch is chosen as the one that gives $w \to 1$ in the limit $\Phi \to 0$, which corresponds to the transition to the initial state of the system $E_{\text{pot}}=0$. Below we show how this general algorithm is implemented in gas dynamics and plasma physics.

\paragraph*{Example 1: Gas Dynamics in a Gravitational Field (Derivation Scheme)}
\label{subsec:gas_derivation}

Consider a stationary flow of an ideal gas in a gravitational field in a pipe of constant cross-section. A detailed thermodynamic analysis, rigorous derivation of this energy balance equation, and its physical interpretation in the context of closed gravitational loops are provided in~\cite{Berezin_Zenodo_Thermo}. The energy conservation law for an elementary layer ($N_s$ particles with mass $m_s$) is written as:
\begin{equation}
\frac{m_s u^2}{2} + \frac{\gamma}{\gamma-1} N_s k_B T + m_s g h = \text{const}.
\end{equation}
Subtracting the state of the layer at the inlet (velocity $u_1$, temperature $T_1$, height $h_1$) and using the continuity equation 
$\rho_1 u_1 = \rho_2 u_2 \Rightarrow u_2 = u_1 w^{-1}$, where $w = \rho_2/\rho_1 = u_1/u_2$, 
as well as the adiabatic relation $T_2 = T_1 w^{\gamma-1}$, we obtain:
\begin{equation}
\frac{m_s u_1^2}{2}\left(1 - \frac{1}{w^2}\right) + 
\frac{\gamma k_B T_1 N_s}{\gamma-1}\left(1 - w^{\gamma-1}\right) + 
m_s g (h_1 - h_2) = 0.
\end{equation}
Introducing the energy variables $\mathcal{K} = \frac{1}{2}m_s u_1^2$, 
$\Theta = \gamma k_B T_1 N_s$, $\Phi = m_s g H = m_s g (h_2 - h_1)$, 
we arrive exactly at the form~\eqref{eq:general_trinomial}. The critical parameter $H_{\max}$ is found from the condition $|x|=R$ for the ultraradical, which physically corresponds to reaching the local speed of sound $u_2 = c_s$.

\paragraph*{Example 2: Plasma Physics (Brief Derivation via Isomorphism)}
\label{subsec:plasma_derivation}

In a stationary ion-acoustic wave in plasma\cite{Dubinov2020}, the energy balance of charged particles in the co-moving system $\xi = z - Vt$ has the same structure~\eqref{eq:general_trinomial}, but with replacement of physical scales:
\begin{equation}
\mathcal{K} = \frac{m_{0} V^2}{2}, \quad 
\Theta = \gamma k_B T_{1}, \quad 
\Phi = Z e (\varphi_2 - \varphi_1).
\end{equation}
Here, instead of the layer velocities $u_1$, $u_2$ relative to the stationary pipe, the particle velocities $v_1, v_2$ relative to the wave velocity $V$ are used.
The continuity equation gives $v = -V/\eta$, where $\eta = n_2/n_1$ is the concentration ratio. The adiabatic state of ions relates $T_2 = T_{1} \eta^{\gamma-1}$. Substitution into~\eqref{eq:general_trinomial} yields:
\begin{equation}
\frac{m_{0} V^2}{2}\left(1 - \frac{1}{\eta^2}\right) + 
\frac{\gamma k_B T_{1}}{\gamma-1}\left(1 - \eta^{\gamma-1}\right) + 
Z e (\varphi_1 - \varphi_2) = 0.
\label{eq:plasma_final}
\end{equation}
The mathematical structure of~\eqref{eq:plasma_final} is identical to the gas dynamics case: the gravitational potential $m_s g H$ is replaced by the electrostatic potential $Z e (\varphi_2 - \varphi_1)$, and the density ratio $\rho_2/\rho_1$ is replaced by the concentration ratio $\eta = n_2/n_1$. Consequently, all analytical results obtained for gas (solution via ultraradical, radius of convergence, critical condition $d\Phi/d\eta=0$) are automatically transferred to plasma.

\paragraph{Remark on the variable adiabatic index and branch selection.}
In real physical processes, the adiabatic index $\gamma$ is not a strict constant: for polyatomic gases, it depends on temperature (due to the excitation of vibrational degrees of freedom), and for water vapor, it even depends on pressure (due to intermolecular interactions near the critical point). As shown in~\cite{Berezin2025_arXiv}, under a smooth variation of the power exponents $\beta_2, \beta_1$, the principal branch of the ultra-radical ($n=0$, as well as $h=0$, $k=0$ for conjugate transformations) remains continuous and analytic. The remaining branches may appear and disappear: for instance, at $\beta_2 = 7$ the equation has seven roots, whereas at $\beta_2 = 5$ it has only five.

As follows from the Sagdeev pseudopotential method~\cite{Dubinov2020}, in problems concerning nonlinear waves (solitons and supersolitons), two positive branches turn out to be physically realizable:
\begin{itemize}
\item \textbf{Light particles (electrons):} their concentration varies according to the $n=0$ branch;
\item \textbf{Heavy particles (ions):} their dynamics is described by the branch $n=-1$.
\end{itemize}

For the trinomial equation to yield two positive real roots, it must possess two positive base points. This is possible exclusively \textbf{beyond the convergence radius of the base expansion} ($|x_1| > R_1$), where the conjugate transformations $j=2$ and $j=3$ come into play, yielding sets of base points of the form:
\[
\left(-\frac{S_2}{S_1}\right)^{\frac{1}{\beta_1 - \beta_2}} \quad \text{and} \quad \left(-\frac{S_1}{S_0}\right)^{\frac{1}{-\beta_1}}.
\]
Thus, the geometric criterion for selecting the conjugate series is not merely a mathematical tool for analytic continuation, but a \textbf{physically necessary mechanism} ensuring the correct description of various plasma components (electrons and ions) during the formation of nonlinear wave structures. This remark illustrates the universality of the proposed formalism: the unified mathematical framework of the ultra-radical allows one to simultaneously and rigorously describe both the thermodynamics of real gases with a variable $\gamma$ and the nonlinear dynamics of multi-component plasmas, where the response of light and heavy particles obeys the laws of different branches of the same trinomial equation.

\section{The Algebra of Master Numbers and the Merge Operation for Master Series}
\label{sec:master_algebra_merge}

A master number of order $\ell$ is a product of factors $(m - a\gamma + b\ell)$. Their count equals $\ell-1$. Therefore, for $\ell<2$, the master number returns 1. Master numbers can be denoted either in the linear form $N(m; a; b; \ell)$ or in the compact form $\mstr{m}{a}{b}{\ell}$, by analogy with the notation of various factorial generalizations.

\begin{equation}
N(m; a; b; \ell) = \mstr{m}{a}{b}{\ell} = \prod_{\gamma = 1}^{\ell-1} (m - a\gamma + b\ell)
\label{eq:master_numbers}
\end{equation}

In expanded form:
\[
\mstr{m}{a}{b}{\ell} = (m - a + b\ell)(m - 2a + b\ell)\cdots\left( m - a(\ell - 1) + b\ell \right).
\]

\subsection{Identities of Master Numbers}
The identities of master numbers are transferred to all generalizations of factorials and powers that are expressed through master numbers.
\[
\mstr{m}{a}{b}{\ell} = (m - a + b\ell)(m - 2a + b\ell)\cdots\left( m - a(\ell - 2) + b\ell \right)\left( m - a(\ell - 1) + b\ell \right)
\]
\[
\mstr{m}{-a}{b-a}{\ell} = \left( m - a(\ell - 1) + b\ell \right)\left( m - a(\ell - 2) + b\ell \right)\cdots(m - 2a + b\ell)(m - a + b\ell)
\]

From this, the symmetry directly follows:
\[
N(m; a; b; \ell) = N(m; -a; b-a; \ell), \quad \mstr{m}{a}{b}{\ell} = \mstr{m}{-a}{b-a}{\ell}.
\]

\textit{Example:}
\[
\mstr{1}{-2}{0}{4} = (1 + 2)(1 + 4)(1 + 6) = 3 \cdot 5 \cdot 7
\]
\[
\mstr{1}{2}{2}{4} = (1 + 8 - 2)(1 + 8 - 4)(1 + 8 - 6) = 7 \cdot 5 \cdot 3
\]

\subsection{Special Cases of Master Numbers}

We show that master numbers \eqref{eq:master_numbers} encompass a wide range of known combinatorial objects \cite{Deza2023, Abramowitz1972, Graham1994}.

\subsubsection*{Factorial}
The ordinary factorial is obtained at $m=1$, $a=1$, $b=1$:
\begin{equation}
\ell! = \mstr{1}{1}{1}{\ell}. \label{eq:factorial}
\end{equation}

\subsubsection*{Multifactorial}
The $a$-fold factorial of a number $n$ is expressed through master numbers as follows:
\begin{equation}
{n!}_{(a)} = \mstr{n-a\ell}{a}{a}{\ell+1}, \qquad \ell = \left\lceil \frac{n}{a} \right\rceil. \label{eq:multifactorial}
\end{equation}

\textit{Examples.} The triple factorial of numbers $5$, $6$, and $7$:
\[
{5!}_{(3)} = \mstr{-1}{3}{3}{3} = 5 \times 2 = 10,
\]
\[
{6!}_{(3)} = \mstr{0}{3}{3}{3} = 6 \times 3 = 18,
\]
\[
{7!}_{(3)} = \mstr{-2}{3}{3}{4} = 7 \times 4 \times 1 = 28.
\]

\subsubsection*{Falling and Rising Factorials}
The falling factorial $(n)_k = n(n-1)\cdots(n-k+1)$ and the rising factorial $n^{(k)} = n(n+1)\cdots(n+k-1)$ are expressed as:
\begin{equation}
(n)_k = \mstr{n+1}{1}{0}{k+1}, \qquad n^{(k)} = \mstr{n-1}{-1}{0}{k+1}. \label{eq:falling_rising}
\end{equation}
The notation and properties of these objects are discussed in detail in \cite{Knuth1992, Knuth1997, Diaz2005}.

\subsubsection*{Binomial Coefficient}
The binomial coefficient through master numbers:
\begin{equation}
\binom{n}{k} = \frac{\mstr{n+1}{1}{0}{k+1}}{k!}. \label{eq:binomial}
\end{equation}
A detailed exposition of the properties of binomial coefficients and their connection with factorials is given in \cite{Knuth2002}.
For example,
\[
\binom{4}{2} = \frac{\mstr{5}{1}{0}{3}}{2!} = \frac{4 \times 3}{2} = 6.
\]

\subsubsection*{Power}
At $a=0$, master numbers take the form:
\begin{equation}
\mstr{m}{0}{b}{\ell} = \prod_{\gamma=1}^{\ell-1} (m + b\ell) = (m + b\ell)^{\ell-1}.
\label{eq:master_a0}
\end{equation}
This form is actively used in generalizations of the Lambert function \cite{Mezo2022, Dubinov2006, Corcino2021}.

\subsubsection*{Generalization}
Thus, all the considered combinatorial objects are special cases of master numbers, which allows using a unified formalism for them.

\subsubsection*{Obtaining Master Numbers Through a Ratio of Gamma Functions}

The master number $N_\ell=\prod_{\gamma=1}^{\ell-1} \left( m - a \gamma + b \ell \right)$ can be expressed in closed form through the gamma function for $a \neq 0$ \cite{Popov2017}:

\begin{equation}
N_\ell = (-a)^{\ell-1} \, 
\frac{\Gamma\!\left( \dfrac{(a-b)\ell - m}{a} \right)}
     {\Gamma\!\left( \dfrac{a - m - \ell b}{a} \right)},
\label{eq:M_L_gamma}
\end{equation}

or in a numerically stable form:
\begin{equation}
N_\ell = a^{\ell-1} \, 
\frac{\Gamma\!\left( \dfrac{m + \ell b}{a} \right)}
     {\Gamma\!\left( \dfrac{m + \ell b}{a} - \ell + 1 \right)}.
\label{eq:M_L_gamma_clean}
\end{equation}

The gamma function has poles at the points $z = 0, -1, -2, \dots$. Consequently, if the argument of the gamma function in the denominator takes a non-positive integer value, then the denominator tends to infinity, and the entire master number becomes zero (provided that the numerator does not tend to infinity at the same time). This gives an explicit analytical criterion for identifying zero terms of the power series.

When using formulas \eqref{eq:M_L_gamma} and \eqref{eq:M_L_gamma_clean}, a situation is possible where the argument of the gamma function simultaneously becomes a pole in both the numerator and the denominator. In this case, the ratio of gamma functions becomes an indeterminate form of type $\infty/\infty$, and the representation of the master number through gamma functions loses meaning. To compute the master number at this point, one should return to its original definition \eqref{eq:master_numbers}.

For the case $a = 0$, master numbers have the form \eqref{eq:master_a0} and cannot be expressed through the gamma function in the manner described above.

\subsection{Integer and Rational Ultraradicals}
\label{subsec:integer_rational}
Although the master-series formalism is universal and applicable to arbitrary complex exponents $a, b \in \mathbb{C}$, the cases where the exponents are rational or integer possess a distinct algebraic and combinatorial structure. We will use the terms \textbf{integer} and \textbf{rational ultraradicals} to describe these cases.

\subsubsection*{Integer Ultraradicals and the Fox--Wright Function}
Consider the particular case of $a = 1$ and an integer $b \in \mathbb{Z}$. In this case, the master numbers \eqref{eq:M_L_gamma_clean} take the form:
\begin{equation}
N_\ell = \frac{\Gamma(m + \ell b)}{\Gamma(m + \ell(b-1) + 1)}.
\end{equation}
Since the arguments of the gamma functions differ by an integer, the ratio of the gamma functions reduces to a falling or rising factorial. The master series for the ultraradical ($m=1$) in this case is written as:
\begin{equation}
M(1; 1; b; x) = 1 + \sum_{\ell=1}^{\infty} \frac{\Gamma(1 + b\ell)}{\Gamma(2 + (b-1)\ell)} \frac{x^\ell}{\ell!}.
\end{equation}
This expression exactly coincides with the \textbf{generalized Fox--Wright hypergeometric function} (or simply the Wright function) ${}_1\Psi_1$:
\begin{equation}
M(1; 1; b; x) = {}_1\Psi_1\left[\begin{array}{c} (1, b) \\ (2, b-1) \end{array} \Big| x\right].
\end{equation}

\subsubsection*{Connection to Classical Combinatorics}
Integer ultraradicals for $m=1$ are closely related to fundamental combinatorial sequences:
\begin{itemize}
\item \textbf{For $b = 1$} (linear equation $y = 1 + xy$): the master numbers are $N_\ell = \ell!$, and the series reduces to a geometric series $M(1; 1; 1; x) = \frac{1}{1-x}$.
\item \textbf{For $b = 2$} (quadratic equation $y = 1 + xy^2$): the series coefficients are proportional to the \textbf{Catalan numbers} $C_\ell = \frac{1}{\ell+1}\binom{2\ell}{\ell}$. The series yields the generating function for the Catalan numbers: $M(1; 1; 2; x) = \frac{1 - \sqrt{1-4x}}{2x}$.
\item \textbf{For $b \ge 3$}: this yields the \textbf{generalized Catalan numbers} (Fuss--Catalan numbers), which describe the combinatorics of $b$-ary trees.
\end{itemize}

\subsubsection*{Rational Ultraradicals}
When the exponents $a$ and $b$ are rational numbers ($a, b \in \mathbb{Q}$), the equation $y^a = 1 + axy^b$ is algebraic. Using the substitution $y = Y^{1/q}$, where $q$ is the common denominator of the fractions $a$ and $b$, it can always be reduced to an equation with integer exponents. Consequently, rational ultraradicals are algebraic functions, and their master series admit cumulative computation.

\subsection{The Merge Operation for Master Numbers}

Master numbers can be merged into composite structures using the operation $\merge$:
\begin{equation}
\mstr{m}{a}{b_1,b_2,\cdots}{\ell_1,\ell_2,\cdots} = \prod_{\gamma = 1}^{\ell_1 + \ell_2 + \cdots - 1}\left( m - a\gamma + b_1\ell_1 + b_2\ell_2 + \cdots \right)
\label{eq:merge_definition}
\end{equation}

\textit{Example.} Merging three master numbers:
\[
\mstr{m}{a}{b_1,b_2,b_3}{\ell_1,\ell_2,\ell_3} = \prod_{\gamma = 1}^{\ell_1 + \ell_2 + \ell_3 - 1}\left( m - a\gamma + b_1\ell_1 + b_2\ell_2 + b_3\ell_3 \right).
\]

For $\ell_1 = \ell_2 = \ell_3 = 1$:
\[
\mstr{m}{a}{b_1,b_2,b_3}{1,1,1} = (m - a + b_1 + b_2 + b_3)(m - 2a + b_1 + b_2 + b_3).
\]

For $m=1$, $a=1$, $b_1=b_2=b_3=1$:
\[
\mstr{1}{1}{1,1,1}{1,1,1} = (1 - 1 + 3)(1 - 2 + 3) = 3 \cdot 2 = 6.
\]

This operation admits pairwise permutations of indices $(b_i, \ell_i)$ and generalizes multidimensional series expansions \cite{Berezin2025_arXiv}.

\subsection{The Merge Operation for Functions}
\label{subsec:merge_functions}

By the \textbf{merge of functions} we mean the merge of their master series: if functions are representable by master series, the merge operation is applied to these series term by term, according to the rule~\eqref{eq:merge_definition}.

If the original functions are compositions of master series (a master series of a master series), then their merge reduces to the merge of the outer master series, while the inner master series remain unchanged. In other words, the merge operation acts at the outer level of the composition and does not affect its inner structure.

Master equations may contain multiple terms with different power exponents $b_i$ and can be solved via master series using a special operation from the class of \emph{nuclear multiplications}, called the \textbf{merge operation for master series}.

\subsubsection*{Main Forms of Master Equations with Multiple Terms}

Consider two main classes of multinomial master equations.

\textbf{Exponential type} ($m=1$, $a=0$):
\begin{equation}
y = e^{x_1 y^{b_1} + x_2 y^{b_2} + \cdots} = M(1; 0; b_1, b_2, \ldots; x_1, x_2, \ldots).
\label{eq:merge_exp}
\end{equation}

\textbf{Power type} ($m=1$, $a \neq 0$):
\begin{equation}
y = \left(1 + a x_1 y^{b_1} + a x_2 y^{b_2} + \cdots\right)^{1/a} = M(1; a; b_1, b_2, \ldots; x_1, x_2, \ldots).
\label{eq:merge_power}
\end{equation}

Their \textbf{logarithmic forms} at $m=0$ are defined analogously:
\begin{equation}
w = x_1 e^{b_1 w} + x_2 e^{b_2 w} + \cdots = M(0; 0; b_1, b_2, \ldots; x_1, x_2, \ldots).
\label{eq:merge_log}
\end{equation}

\textbf{Example.} The equation $w = x \cosh(bw)$ can be rewritten as:
\begin{equation}
w = \frac{x \cdot e^{bw} + x \cdot e^{-bw}}{2} = M\!\left(0; 0; b, -b; \frac{x}{2}, \frac{x}{2}\right).
\label{eq:merge_cos_example}
\end{equation}

\subsubsection*{Reduction to a Single-Variable Series}

Multiplying series with multiple independent variables $x_1, x_2, \ldots$ is inconvenient for analysis. To obtain a familiar power series in a single variable, the substitution is applied:
\begin{equation}
x_1 = z_1 Z, \quad x_2 = z_2 Z, \quad \ldots
\label{eq:merge_substitution}
\end{equation}
where $z_i$ are fixed coefficients and $Z$ is the unified expansion variable.

\subsubsection*{Formal Definition of the Merge Operation}

The merge operation for master series is viewed as the multiplication of several power series, in which the master numbers themselves are \emph{not multiplied} but \emph{merged} according to the rule~\eqref{eq:merge_definition}:
\begin{equation}
N(m; a; b_1; \ell_1) \merge N(m; a; b_2; \ell_2) = N(m; a; b_1, b_2; \ell_1, \ell_2).
\label{eq:merge_numbers_rule}
\end{equation}

In the merging of series, the zeroth term is always equal to $1$, regardless of the value of the parameter $m$. In the case $m=0$ (logarithmic form), after completing the merge operation, $1$ must be subtracted from the result.

\subsubsection*{Explicit Form of Merging Two Series}

The general formula for merging several master series is
\begin{equation}
M(m; a; b_1, b_2, \ldots; x_1, x_2, \ldots) = M(m; a; b_1; x_1) \merge M(m; a; b_2; x_2) \merge \cdots
\label{eq:merge_general}
\end{equation}
Let us illustrate the merge operation with an example of two master series:
\begin{align}
&\left(1 + \frac{N_{1,1} z_1 Z}{1!} + \frac{N_{1,2} z_1^2 Z^2}{2!} + \frac{N_{1,3} z_1^3 Z^3}{3!} + \cdots\right) \merge \notag \\
&\qquad \merge \left(1 + \frac{N_{2,1} z_2 Z}{1!} + \frac{N_{2,2} z_2^2 Z^2}{2!} + \frac{N_{2,3} z_2^3 Z^3}{3!} + \cdots\right) = \notag \\
&= 1 + \frac{(N_{1,1} z_1 + N_{2,1} z_2) Z}{1!} + \notag \\
&\quad + \frac{\left(N_{1,1} \merge N_{2,1} \cdot z_1 z_2 \cdot 2 + N_{1,2} z_1^2 + N_{2,2} z_2^2\right) Z^2}{2!} + \notag \\
&\quad + \frac{\left(N_{1,1} \merge N_{2,2} \cdot z_1 z_2^2 \cdot 3 + N_{1,2} \merge N_{2,1} \cdot z_1^2 z_2 \cdot 3 + N_{1,3} z_1^3 + N_{2,3} z_2^3\right) Z^3}{3!} + \cdots
\label{eq:merge_explicit}
\end{align}

\subsubsection*{Compact merge formula}
The explicit expression~\eqref{eq:merge_explicit} generalizes naturally to the merge of an arbitrary number $k$ of master series. The coefficient $\eta_\ell$ of $Z^\ell$ in the expansion
\begin{equation}
M(m; a; b_1, \ldots, b_k; z_1 Z, \ldots, z_k Z) = m + \sum_{\ell=1}^{\infty} \frac{\eta_\ell}{\ell!} Z^\ell
\label{eq:merge_power_series}
\end{equation}
is given by the compact formula:
\begin{equation}
\eta_\ell = \sum_{\substack{r_1 + \cdots + r_k = \ell \\ r_i \geq 0}} \frac{\ell!}{r_1! \cdots r_k!} z_1^{r_1} \cdots z_k^{r_k} \prod_{\gamma=1}^{\ell-1} \left(m - a\gamma + \sum_{i=1}^k b_i r_i\right),
\label{eq:eta_compact}
\end{equation}
where the empty product (for $\ell=0,1$) is defined as unity. If the common factor $Z$ of the independent arguments is omitted, it is assumed to be $1$, and formula~\eqref{eq:merge_power_series} reduces to the standard notation of the master series with arguments $x_i = z_i$.

\subsubsection*{Super-master merging}
For super-master series with a common parameter $c$, the merge operation yields:
\begin{equation}
S(m; a; b_1, \ldots, b_k; z_1 Z, \ldots, z_k Z; c) = m + c\sum_{\ell=1}^{\infty} \frac{\eta_\ell^{(c)}}{\ell!} Z^\ell,
\label{eq:super_merge_series}
\end{equation}
where the coefficients $\eta_\ell^{(c)}$ are given by the formula:
\begin{equation}
\eta_\ell^{(c)} = \sum_{\substack{r_1 + \cdots + r_k = \ell \\ r_i \geq 0}} \frac{\ell!}{r_1! \cdots r_k!} z_1^{r_1} \cdots z_k^{r_k} \prod_{\gamma=1}^{\ell-1} \left(cm - a\gamma + \sum_{i=1}^k b_i r_i\right),
\label{eq:super_eta_compact}
\end{equation}
where the empty product (for $\ell=0,1$) is set equal to unity. At $c=1$, formulas~\eqref{eq:super_merge_series}--\eqref{eq:super_eta_compact} reduce to~\eqref{eq:merge_power_series}--\eqref{eq:eta_compact}.

Substituting expression~\eqref{eq:super_eta_compact} into the expansion~\eqref{eq:super_merge_series} causes the factorial $\ell!$ in the numerator of the multinomial coefficient to cancel with the $\ell!$ in the denominator of the series. Consequently, the $\ell$-th term of the super-master series assumes a compact form:
\begin{equation}
C_\ell=c Z^\ell \sum_{\substack{r_1 + \cdots + r_k = \ell \\ r_i \geq 0}} \frac{z_1^{r_1}}{r_1!} \cdots \frac{z_k^{r_k}}{r_k!} \prod_{\gamma=1}^{\ell-1} \left(cm - a\gamma + \sum_{i=1}^k b_i r_i\right).
\label{eq:super_merge_simplified}
\end{equation}

The algorithmic computation of this term is formalized in Algorithm~\ref{alg:CL}. To avoid combinatorial overflow and ensure $O(1)$ amortized time per combination, the algorithm utilizes an iterative ``odometer'' mechanism for generating weak compositions.

\begin{algorithm}[H]
\caption{Computation of the $L$-th term of the merged super-master series (CL)}
\label{alg:CL}
\begin{algorithmic}[1]
\Require $L \in \mathbb{N}$, $k \in \mathbb{N}$, $a, c, m \in \mathbb{C}$; arrays $b, x$ of length $k$
\Ensure The value of the $L$-th term (without the $Z^L$ factor)
\State $S \gets 0$, $H \gets 0$
\State Initialize array $r$ of length $k$ with zeros
\While{true}
    \While{$S \le L$}
        \State $r[k] \gets L - S$, $P \gets 1$, $B \gets 0$
        \For{$g = 1$ \textbf{to} $k$}
            \State $P \gets P \cdot \frac{x[g]^{r[g]}}{r[g]!}$
            \State $B \gets B + b[g] \cdot r[g]$
        \EndFor
        \For{$\gamma = 1$ \textbf{to} $L-1$} \Comment{Master number product}
            \State $P \gets P \cdot (cm - a\gamma + B)$
        \EndFor
        \State $H \gets H + P$
        \State $r[k-1] \gets r[k-1] + 1$
        \State $S \gets S + 1$
    \EndWhile
    \State $h \gets k - 1$
    \While{$h > 1$ \textbf{and} $r[h] = 0$}
        \State $h \gets h - 1$
    \EndWhile
    \If{$h = 1$}
        \State \Return $c \cdot H$
    \EndIf
    \State $r[h] \gets 0$
    \State $r[h-1] \gets r[h-1] + 1$
    \State $S \gets 0$
    \For{$i = 1$ \textbf{to} $k-1$}
        \State $S \gets S + r[i]$
    \EndFor
\EndWhile
\end{algorithmic}
\end{algorithm}

Thus, the merged super-master series itself takes a remarkably compact form:
\begin{equation}
S(m; a; b_1, \ldots, b_k; z_1 Z, \ldots, z_k Z; c) = m + \sum_{\ell=1}^{\infty} C_\ell.
\label{eq:super_merge_final}
\end{equation}

As a substantive illustration of the applicability of the proposed formalism, we consider a problem that has historically served as a serious challenge for analytical methods: the solution of a quintic polynomial equation. In terms of the merge operation for super-master series, the solution of an equation of the form
\begin{equation}
y^5 = 1 + 5 \cdot 0.01 y^4 + 5 \cdot 0.001 y^3 + 5 \cdot 0.005 y^2 + 5 \cdot 0.002 y
\label{eq:quintic_example}
\end{equation}
can be written in a compact form that structurally unifies the classical power series of Lambert and Euler for trinomial equations with their modern rigorous generalization to arbitrary powers via Bell polynomials and the Fox--Wright function \cite{Belkic2019}:
\begin{equation}
Y_n^{c} = e^{c f_n} \cdot \text{MergeSuperMaster}(1, 5, \mathbf{b}, \mathbf{x}, c, M_{\max})[1],
\label{eq:quintic_solution}
\end{equation}
where the parameter vectors are given by $\mathbf{b} = [4, 3, 2, 1]$ and $\mathbf{x} = [0.01 e^{4 f_n}, 0.001 e^{3 f_n}, 0.005 e^{2 f_n}, 0.002 e^{f_n}]$, and $f_n$ is defined by formula~\eqref{eq:f_def}.

This example demonstrates how a unified algorithmic approach allows one to constructively work with equations that previously required cumbersome specific substitutions or numerical methods. In physical applications, it is critically important to unambiguously specify which particular solution branch is physically realizable. As shown in the introduction of the ``root passport'', the proposed formalism allows this to be done even for equations with irrational power exponents. Moreover, for trinomial equations, we have developed a geometric criterion that ensures continuous analytic continuation of branches beyond the radius of convergence.

\subsubsection*{Multidimensional Convergence Domain of the Merged Series}
The merged super-master series $y = S(1; a; b_1, \ldots, b_k; x_1, \ldots, x_k; 1)$ is the solution to the multinomial equation
\begin{equation}
y^a = 1 + a \sum_{i=1}^k x_i y^{b_i}.
\label{eq:multinomial_eq}
\end{equation}
According to the theory of analytic functions, the radius of convergence of a power series representing a root of an algebraic (or generalized algebraic) equation $F(y, \mathbf{x}) = 0$ is bounded by the nearest singular point, where the implicit function theorem fails, i.e., where $\frac{\partial F}{\partial y} = 0$.
For equation~\eqref{eq:multinomial_eq}, this condition takes the form:
\begin{equation}
a y^{a-1} - a \sum_{i=1}^k b_i x_i y^{b_i - 1} = 0 \quad \Longrightarrow \quad y^a = \sum_{i=1}^k b_i x_i y^{b_i}.
\label{eq:singularity_condition}
\end{equation}
Combining~\eqref{eq:multinomial_eq} and~\eqref{eq:singularity_condition}, we obtain a system that defines the discriminant manifold (the boundary of the convergence domain) in the $k$-dimensional space of arguments:
\begin{equation}
\begin{cases}
\displaystyle \sum_{i=1}^k x_i y^{b_i} = \frac{y^a - 1}{a}, \\[10pt]
\displaystyle \sum_{i=1}^k b_i x_i y^{b_i} = y^a.
\end{cases}
\label{eq:discriminant_system}
\end{equation}

\textbf{Convergence along a given ray.} If the arguments vary proportionally along a ray $x_i = z_i Z$ (where $z_i$ are fixed directions and $Z$ is a scalar parameter), system~\eqref{eq:discriminant_system} reduces to a single equation for the critical value $y_c$:
\begin{equation}
y_c^a \left( 1 - a \frac{\sum_{i=1}^k z_i y_c^{b_i}}{\sum_{i=1}^k b_i z_i y_c^{b_i}} \right) = 1.
\label{eq:yc_equation}
\end{equation}
Finding $y_c$, we obtain the exact radius of convergence $R = |Z_c|$ along this ray:
\begin{equation}
Z_c = \frac{y_c^a}{\sum_{i=1}^k b_i z_i y_c^{b_i}}.
\label{eq:Zc_formula}
\end{equation}
\textit{Check for $k=1$.} For $k=1$ ($z_1=1, b_1=b$), equation~\eqref{eq:yc_equation} yields $y_c^a(1 - a/b) = 1$, whence $y_c = (1 - a/b)^{-1/a}$. Substitution into~\eqref{eq:Zc_formula} gives $Z_c = \frac{1}{b} y_c^{a-b} = \frac{1}{b}\left(\frac{b-a}{b}\right)^{\frac{b-a}{a}}$, which exactly coincides with the classical radius of convergence $R = \frac{|1 - a/b|^{b/a}}{|b - a|}$ derived earlier for the single master series.

\textbf{Convergence boundary for $k=2$.} In the case of merging two series ($k=2$), system~\eqref{eq:discriminant_system} is uniquely solvable for $x_1$ and $x_2$. The boundary of the two-dimensional convergence domain is given by the parametric curve:
\begin{equation}
x_1(y) = \frac{y^a(a - b_2) + b_2}{a(b_1 - b_2)} y^{-b_1}, \qquad x_2(y) = \frac{y^a(b_1 - a) - b_1}{a(b_1 - b_2)} y^{-b_2}.
\label{eq:k2_boundary}
\end{equation}
This result strictly describes the geometry of the convergence domain for merged series and demonstrates how the combinatorial structure of master numbers adapts to the multidimensional space of arguments.

\paragraph{Remark on practical convergence testing.}
While the discriminant system \eqref{eq:discriminant_system} provides an exact characterization of the convergence boundary, it is not directly tractable in practice: for $k \ge 3$, solving this system is generally no simpler than solving the original multinomial equation itself. However, for practical branch identification, it often suffices to determine which of the canonical transformations yields a convergent series along a given physical trajectory. This motivates the search for practically computable, iterative criteria.

A natural heuristic emerges from the geometry of multidimensional Reinhardt domains. At the origin ($x_i = 0$), the roots are simply the $a$-th roots of the base ratio. As one argument, say $x_1$, increases towards its individual one-dimensional radius $R_1$, the series remains convergent. If a second argument $x_2$ is then turned on and increased, the maximum allowable modulus $|x_1|_{\max}$ must be decreased accordingly to preserve convergence. We conjecture that an explicit, monotonic trade-off relation between $|x_1|$ and $|x_2|$ exists, allowing one to compute the effective radius of one argument while the others are held fixed. Rigorous derivation of such scaling laws and the development of an iterative boundary-tracing algorithm for $k \ge 3$ remain open problems and will be the subject of our future work.

\subsubsection*{Criterion for Asymptotic Term Exclusion and the Newton Polygon}
For equations with real exponents, the convergence question is global for the entire canonical transformation: the series either converges for all branches or diverges for all of them. This property allows one to discard inefficient transformations at the preliminary analysis stage based on the term exclusion criterion, thereby avoiding computationally expensive branch-by-branch enumeration.

Consider the general form of a multinomial equation:
\begin{equation}
\sum_{j=0}^{k} S_j Y^{\beta_j} = 0, \qquad \beta_k > \beta_{k-1} > \dots > \beta_1 > \beta_0 \ge 0.
\label{eq:multinomial_general}
\end{equation}
The extreme terms ($j=k$ and $j=0$) always remain leading and are not subject to exclusion. In the case of a trinomial ($k=2$), the question reduces to the analysis of a single middle term ($j=1$). Let us introduce the differences of the power exponents: $\alpha_2 = |\beta_1 - \beta_0|$, $\alpha_1 = |\beta_2 - \beta_0|$, $\alpha_0 = |\beta_2 - \beta_1|$, as well as the negligibility indicators $T_j = |S_j / \alpha_j|^{\alpha_j}$. The condition that the middle term $S_1 Y^{\beta_1}$ is asymptotically negligible compared to the extreme terms reduces to the balance criterion:
\begin{equation}
T_1 < T_2 T_0.
\label{eq:dominance_trinomial}
\end{equation}
If this condition is satisfied, the middle term is excluded from consideration, and all roots are determined via a single binary transformation, whose base points form a regular $\alpha_1$-gon with vertices at $(-S_0/S_2)^{1/\alpha_1}$.

This approach naturally generalizes to an arbitrary number of terms and is intuitively interpreted in terms of the \textbf{Newton polygon}. Each pair of adjacent (non-excluded) terms $S_i Y^{\beta_i}$ and $S_j Y^{\beta_j}$ forms its own cluster of roots, asymptotically located near the vertices of a regular $|\beta_i - \beta_j|$-gon.

\textbf{Example.} Consider the equation $4Y^{12} + 3Y^7 + 2Y^3 + 5 = 0$. Due to the smallness of the coefficients at the middle terms, they are subject to exclusion when forming the Newton polygon. As a result, all 12 roots are located near the vertices of a regular dodecagon $(-5/4)^{1/12}$. 
If the coefficient of $Y^7$ is increased to 30, considering the equation $4Y^{12} + 30Y^7 + 2Y^3 + 5 = 0$, the exclusion condition for this term is violated. The set of roots then splits into two clusters: 5 roots group near the vertices of a pentagon (from the pair $4Y^{12}$ and $30Y^7$, since $12-7=5$) and 7 roots --- near the vertices of a heptagon (from the pair $30Y^7$ and $5$, since $7-0=7$). 
Similarly, increasing the coefficient of $Y^3$ to 20 in the equation $4Y^{12} + 3Y^7 + 20Y^3 + 5 = 0$ leads to the splitting of the dodecagon into triangles and nonagons. 
In the case where none of the terms can be excluded, for example, in the equation $4Y^{12} + 30Y^7 + 32Y^3 + 5 = 0$, each pair of consecutive terms forms its own polygon: the pair $4Y^{12}$ and $30Y^7$ forms a pentagon, the pair $30Y^7$ and $32Y^3$ forms a square, and the pair $32Y^3$ and $5$ forms an equilateral triangle.

It should be emphasized that there is a fundamental connection between this geometric representation and the analytical properties of the series: as the arguments move away from the boundary of the radius of convergence, the roots of the equation asymptotically approach their base points, forming the regular polygonal structures described above.

\subsubsection*{Geometry of base points}
\label{subsec:base_points_geometry}

\textbf{Base points} are defined as the values of the solution at $x=0$, that is, the roots of the degenerate equations obtained from \eqref{eq:canonical_pair} at $x=0$. The parameter $m$ in the master series $M(m; a; b; x)$ has a deep geometric meaning: it determines the arrangement of the solution branches in the neighborhood of $x = 0$.

\textbf{Case $m=1$ (ultraradical).} As $x \to 0$, the equation $y^a = 1 + axy^b$ reduces to $y^a = 1$, and all branches tend to the roots of unity:
\begin{equation}
y_n \to e^{2\pi i n/a}, \qquad n \in \mathbb{Z}.
\end{equation}
For \emph{integer} $a$, these points form a regular $a$-gon on the unit circle $|y| = 1$. For non-integer real $a$, the set $\{e^{2\pi i n/a}\}_{n\in\mathbb{Z}}$ is finite (if $a$ is rational) or everywhere dense on the circle (if $a$ is irrational); for complex $a$, the points lie on the circle $|y|=1$, but do not fill it uniformly.

\textbf{Case $m=0$ (ultralogarithm).} As $x \to 0$, the equation $e^{aw} = 1 + axe^{bw}$ reduces to $e^{aw} = 1$, and all branches tend to the points:
\begin{equation}
w_n \to \frac{2\pi i n}{a}, \qquad n \in \mathbb{Z}, \qquad a \neq 0.
\end{equation}
For \emph{real} $a \neq 0$, these points lie on the imaginary axis $\operatorname{Re}(w) = 0$, forming a uniform lattice with step $2\pi/a$. For complex $a$, they lie on the line $\{2\pi i t/a : t \in \mathbb{R}\}$, passing through the origin at an angle $\pi/2 - \arg(a)$.

\textbf{Exponential mapping.} The fundamental identity $y = e^w$ establishes a geometric correspondence between these two pictures: the exponential maps the imaginary axis (where the base points of the ultralogarithm lie) onto the unit circle (where the base points of the ultraradical lie). Thus, the base points of the ultralogarithm $w_n = 2\pi i n/a$ are exactly the preimages of the base points of the ultraradical $y_n = e^{2\pi i n/a}$ under the mapping $y = e^w$.

This geometric correspondence reflects the algebraic structure of the recurrence for derivatives: the polynomials $P_\ell(\xi; 1)$ and $P_\ell(\xi; 0)$ are related through the Stirling numbers of the second kind (see~\eqref{eq:stirling_relation}), and this connection is an algebraic reflection of the transition from power branching to logarithmic branching.

\textbf{Behavior in the neighborhood of base points.} In the neighborhood of $x=0$, the solutions of the ultraradical and ultralogarithm have the form
\begin{equation}
y_n(x) = y_n(0) \cdot M(1; a; b; x), \qquad w_n(x) = w_n(0) + M(0; a; b; x),
\end{equation}
where $y_n(0) = e^{2\pi i n/a}$, $w_n(0) = 2\pi i n/a$, and $M(m;a;b;x)$ is the corresponding master series \eqref{eq:maclaurin_master}. In other words, in the ultraradical the base point is \emph{multiplied} by the master series, while in the ultralogarithm it is \emph{added} to it. This is a direct consequence of the identity $y = e^w$: the exponential turns addition into multiplication.

\subsubsection*{Analytic Continuation}

Beyond the radius of convergence, similarly to the trinomial case, transformations of the original equations are applied. However, for equations \textbf{with four or more terms}, the number of possible canonical transformations grows significantly and exceeds three, which requires a generalization of the geometric criterion for selecting the conjugate series.

\subsubsection*{Other Operations of the Class of Nuclear Multiplications}
In the preceding sections, we introduced the term \emph{nuclear multiplications} to describe operations on master series. However, it requires clarification why we use the specific term \emph{merge operation} rather than simply ``a generalization of the super-master to multiple terms.'' The point is that when combining power series, their coefficients --- the master numbers --- are not multiplied as ordinary scalars. Rather, their internal structure interacts --- their ``kernel'' (a finite product of linear factors). Depending on the algebraic structure of the original equation, the laws governing the interaction of these kernels may differ. Thus, the ``merge'' is merely one, albeit fundamental, specific law of master-number interaction, generated by equations of the form $y^a = 1 + a \sum x_i y^{b_i}$. We propose to call the entire class of possible laws of master-number interaction under series combination \emph{nuclear multiplications}.

A natural question arises: do there exist other laws of nuclear multiplication distinct from the merge? It turns out that not all types of polynomial and exponential-power equations reduce to the standard merge. The structure of the equation dictates its own unique law of kernel interaction.

Consider, for example, an equation containing the square of a sum in the exponent or argument:
\begin{equation}
y^{x_1 y^{-b_1} + x_2 y^{-b_2} + \dots} = e^{g(x_1 + x_2 + \dots)^2},
\label{eq:non_merge_example}
\end{equation}
where $g$ is a parameter. In such equations, the combinatorics of the coefficients in the series expansion changes: instead of the linear additivity of indices characteristic of the merge, a quadratic dependence and cross-terms of the form $x_i x_j$ arise.

At present, it remains an open question whether every equation of this type can be reduced by algebraic substitutions to a canonical form solvable via the standard merge of master series. If this is not possible, such equations define a new, independent class of problems. Their solution will require investigating and formalizing the specific law of nuclear multiplication that operates on master series within these structures, establishing new rules for combining their kernels. The investigation of this class of operations constitutes an important direction for the further development of the proposed formalism.

\section{Convergence rate of series and hybrid algorithm with Riemann sheet anchoring}
\label{sec:practical_convergence_hybrid}

A well-known practical limitation of any power series, including master series, is the slowdown of convergence as the argument approaches the boundary of the circle of convergence ($|x| \to R$). In theory, such series are ideal for analytical derivations; however, in practice, this creates serious computational difficulties. 

A vivid illustration of this is provided by the quintic equation considered above~\eqref{eq:quintic_example}. Numerical experiments in Maple showed that at standard precision \texttt{Digits := 10}, all five complex roots converge in just 5--6 terms of the series (the \texttt{MergeSuperMaster} output fixes $L=5$ or $L=6$). This occurs because in this example, the moduli of the arguments $|x_j|$ are significantly smaller than the radius of convergence $R$. However, if the argument is close to $R$, or if high precision is required (50--100 significant digits), the number of required series terms increases to hundreds or even thousands. Even when using optimized code, direct summation of such a number of terms becomes computationally expensive, which is absolutely unacceptable when plotting in physical problems that require calculating roots for thousands of different argument values (for example, when calculating potential profiles or phase portraits).

To overcome this limitation, we have developed a \textbf{hybrid algorithm} that combines rapid series summation and iterative refinement. The algorithm works as follows: first, the super-master series is computed with a strict limit on the maximum number of terms $M_{\max}$. If the series reaches machine precision in $L < M_{\max}$ steps, the process terminates. However, if $L \ge M_{\max}$ (which implies forced truncation of the series and no guarantee of accuracy), the algorithm automatically switches to iterative refinement mode, using the obtained truncated value as the initial approximation $y_0$.

Here arises a key problem, which is solved by the verification index $u$. When iteratively solving equations with fractional or complex power exponents (such as $y^a = 1 + axy^b$), the desired root may be verified in the original equation not through the principal branch of the complex logarithm. Consequently, if the iterative formulas (Newton's, Schröder's) perform exponentiation to fractional or complex powers via the principal branch, the identity will simply not be achieved: the iterations will either converge to a different root or diverge. To correctly find the root, it is necessary to explicitly specify in the iterative formulas the branch of exponentiation (the index $u$) that leads to the stated identity.

To guarantee the preservation of the branch topology, we introduce a \textbf{generalized power function with rigid fixation of the Riemann sheet}:
\begin{equation}
\text{pow}(y, p; u) = \exp\!\Big( p \big[ \ln|y| + i(\arg y + 2\pi u) \big] \Big),
\label{eq:pow_anchored}
\end{equation}
where $u$ is the verification index obtained from the ``root passport'' at the stage of the super-master series or the geometric criterion. 

When constructing the iterative process (in this case, an adaptive Newton--Schröder method), all nonlinear terms of the equation $F(y) = 0$ and its derivatives $F'(y), F''(y)$ are computed exclusively through the function~\eqref{eq:pow_anchored}. Thus, the parameter $u$ acts as a \textbf{topological anchor}, which forcibly holds the iterative trajectory on the exact same sheet of the Riemann surface where the initial approximation was found. The iterative process takes the form:
\begin{equation}
y_{j+1} = y_j - \frac{F(y_j) F'(y_j)}{F'(y_j)^2 - F(y_j) F''(y_j)},
\label{eq:schroeder_step}
\end{equation}
where the choice between the classical Newton method (when $|F'(y_j)| > \varepsilon$) and the Schröder method (for small derivatives) ensures the stability of convergence.

The proposed hybrid algorithm solves the problem of slow series convergence in a radical way. The super-master series acts as a fast ``navigator'', delivering the initial approximation to the required region of the complex plane in just a few terms and unambiguously determining the required Riemann sheet via the index $u$. Next, the iterative refiner comes into play, which, thanks to the ``anchoring'' via~\eqref{eq:pow_anchored}, reaches machine precision in 3--5 steps without risking slipping onto a foreign branch.

\section{Non-limit cases of the master series and classical functions}
\label{sec:classical_functions}

\subsection{Connection with the Inverse Hyperbolic Sine.} 
The defining equation of the ultralogarithm for $a=2$ and $b=1$ takes the form $e^{2w} - 2xe^w - 1 = 0$. 
Dividing both sides by $2e^w$ directly yields the canonical definition of the hyperbolic sine:
\begin{equation}
\frac{e^w - e^{-w}}{2} = x \quad \Longrightarrow \quad \sinh(w) = x.
\end{equation}
Consequently, the ultralogarithm $M(0; 2; 1; x)$ exactly represents the function $w = \operatorname{arsinh}(x)$. 

This demonstrates a fundamental advantage of the ultralogarithm formalism over the classical definition $\operatorname{arsinh}(x) = \ln(x + \sqrt{x^2+1})$. While the classical expression requires the composition of multivalued functions (logarithm and square root) and the use of a discrete $\pm$ sign to describe all branches, the equation $e^{2w} - 2xe^w - 1 = 0$ constitutes a single, indivisible analytic object. The branch parameter $n \in \mathbb{Z}$ in the master series $M(0; 2; 1; x)$ unifies the choice of the $\pm$ sign and the shift by $2\pi i k$ into a single topological index, ensuring strict and continuous tracking of all sheets of the Riemann surface for $\operatorname{arsinh}(x)$.

\subsection{Connection with hyperbolic sine, cosine, and arctangent.}
The master series formalism naturally encompasses all classical hyperbolic functions. Consider the master series $M(1;0;0;x)$, which at $a=0, b=0$ represents the exponential function:
\begin{equation}
M(1;0;0;x) = e^x = 1 + x + \frac{x^2}{2!} + \frac{x^3}{3!} + \frac{x^4}{4!} + \cdots
\end{equation}
Separating this series into even and odd terms yields exactly the hyperbolic cosine and sine:
\begin{align}
\cosh(x) &= 1 + \frac{x^2}{2!} + \frac{x^4}{4!} + \cdots = \frac{M(1;0;0;x) + M(1;0;0;-x)}{2}, \\
\sinh(x) &= x + \frac{x^3}{3!} + \frac{x^5}{5!} + \cdots = \frac{M(1;0;0;x) - M(1;0;0;-x)}{2}.
\end{align}

An analogous decomposition holds for the logarithm $M(0;1;0;x) = \ln(1+x)$:
\begin{equation}
\ln(1+x) = x - \frac{x^2}{2} + \frac{x^3}{3} - \frac{x^4}{4} + \frac{x^5}{5} - \cdots
\end{equation}
The odd terms of this series form the inverse hyperbolic tangent:
\begin{equation}
\operatorname{artanh}(x) = x + \frac{x^3}{3} + \frac{x^5}{5} + \cdots = \frac{M(0;1;0;x) - M(0;1;0;-x)}{2}.
\end{equation}

\subsection{Trigonometric and Hyperbolic Functions as Compositions of Master Series}
The master series formalism encompasses not only hyperbolic but also trigonometric functions through composition with basic master series. The key observation: trigonometric functions are expressed through the binomial root $M(1; -1; 0; z) = (1-z)^{-1}$:
\begin{align}
\sec(x) &= \frac{1}{\cos(x)} = M(1; -1; 0; 1-\cos(x)), \\
\csc(x) &= \frac{1}{\sin(x)} = M(1; -1; 0; 1-\sin(x)).
\end{align}
Other trigonometric ratios are constructed through these functions:
\begin{align}
\tan(x) &= \sin(x) \cdot \sec(x) = \sin(x) \cdot M(1; -1; 0; 1-\cos(x)), \\
\cot(x) &= \cos(x) \cdot \csc(x) = \cos(x) \cdot M(1; -1; 0; 1-\sin(x)).
\end{align}
Similarly, inverse trigonometric functions are expressed through the logarithm:
\begin{equation}
\arctan(x) = \frac{1}{2i} \ln\!\left(\frac{1+ix}{1-ix}\right) = \frac{1}{2i} \left[ M(0; 1; 0; ix) - M(0; 1; 0; -ix) \right].
\end{equation}
Thus, the master formalism provides a unified analytical apparatus for representing all classical trigonometry, where each function is expressed through the composition of basic master series.

\section{Derivatives of ultrafunctions of arbitrary order}
\label{sec:derivatives}
Consider the canonical equations of the ultra-radical and ultra-logarithm:
\begin{equation}
y^a = 1 + a x y^b, \qquad e^{aw} = 1 + a x e^{bw},
\label{eq:canonical_pair}
\end{equation}
where \(a, b, x \in \mathbb{C}\). The solutions of these equations are the branches \(y_n(x)\) and \(w_n(x)\), \(n \in \mathbb{Z}\), related by the identity \(y_n = e^{w_n}\).

When \(b = 0\), equations \eqref{eq:canonical_pair} reduce to \(y^a = 1 + ax\) and \(e^{aw} = 1 + ax\), so the ultra-radical and ultra-logarithm generalize the classical power \((1+ax)^{1/a}\) and logarithmic \((1/a)\ln(1+ax)\) functions. The prefix <<ultra>> emphasizes that we transition from binomial equations (root, logarithm) to trinomial ones; the limiting case \(b = 0\) demonstrates how this transition reverses, returning us to the classical functions.

Let us introduce the \textbf{derivative of the ultralogarithm} (which also serves as the logarithmic derivative of the ultraradical) $\psi$, and a unified \textbf{dimensionless variable} $\xi$:
\begin{equation}
\psi = \frac{y'}{y} = w', \qquad \xi = (a-b)\,b\,x\,\psi.
\label{eq:psi_xi_def}
\end{equation}
Also define auxiliary parameters depending only on the power exponents:
\begin{equation}
\sigma = 2b - a, \qquad \tau = (a-b)b.
\label{eq:sigma_tau_def}
\end{equation}
These two numbers are the coefficients of the quadratic equation
\begin{equation}
\zeta^2 = \tau + \sigma \zeta,
\label{eq:master_equation_exponents}
\end{equation}
whose roots are exactly \(\zeta_1 = b\) and \(\zeta_2 = b - a\). The quadratic trinomial \(\tau + \sigma\xi - \xi^2\), arising in the partial differential equation \eqref{eq:pde_g} for the generating function, factors through these roots:
\[
\tau + \sigma\xi - \xi^2 = (b - \xi)(\xi - b + a) = -(\xi - \zeta_1)(\xi - \zeta_2).
\]
This factorization allows one to decompose the integrand in the characteristic equation
\[
\frac{d\xi}{\tau + \sigma\xi - \xi^2}
\]
into partial fractions and integrate by the method of characteristics in elementary functions (see~\eqref{eq:characteristics}).
\subsection{First derivatives}
Differentiating the canonical equations with respect to $x$, we find the first derivatives:
\begin{equation}
y' = \frac{y}{y^{a-b} - b x}, \qquad w' = \frac{e^{bw}}{1 + (a-b)x e^{bw}}.
\label{eq:first_derivatives}
\end{equation}
\subsection{Unified derivative polynomial}
The key object describing the higher derivatives of both equations is the family of polynomials $P_\ell(\xi; m)$ of degree $\ell$ in the variable $\xi$, parameterized by the value of $m$.

\textbf{Definition.} The polynomials $P_\ell(\xi; m)$ are defined by the initial condition $P_0(\xi; m) = 1$ and the recurrence relation:
\begin{equation}
\boxed{P_{\ell+1}(\xi; m) = \bigl[m + (\ell+1)\sigma - (\ell+1)\xi\bigr] P_\ell(\xi; m) + P_\ell'(\xi; m) \bigl[\tau + \sigma \xi - \xi^2\bigr],}
\label{eq:P_unified_recur}
\end{equation}
where $P_\ell'(\xi; m)$ denotes the derivative of the polynomial with respect to the variable $\xi$.
\subsection{Theorems on higher derivatives}
Having at our disposal the unified polynomial $P_\ell(\xi; m)$, the formulas for higher derivatives take an extremely symmetric form.

\textbf{Theorem 1 (Ultra-radical).} The derivatives of the solution $y(x)$ are expressed through the polynomial $P_\ell(\xi; 1)$:
\begin{equation}
\boxed{y^{(\ell+1)} = y' \, \psi^\ell \, P_\ell(\xi; 1), \qquad \ell \geq 0.}
\label{eq:y_n_theorem}
\end{equation}

\textbf{Theorem 2 (Ultra-logarithm).} The derivatives of the solution $w(x)$ are expressed through the polynomial $P_\ell(\xi; 0)$:
\begin{equation}
\boxed{w^{(\ell+1)} = \psi^{\ell+1} \, P_\ell(\xi; 0), \qquad \ell \geq 0.}
\label{eq:w_n_theorem}
\end{equation}
\textbf{Proof.} We differentiate the identity \eqref{eq:y_n_theorem} with respect to $x$:
\begin{equation}
y^{(\ell+2)} = \frac{d}{dx}\bigl[y' \psi^\ell P_\ell(\xi; 1)\bigr] = y'' \psi^\ell P_\ell + y' \cdot \ell \psi^{\ell-1} \psi' P_\ell + y' \psi^\ell P_\ell'(\xi; 1) \xi'.
\end{equation}
From equation \eqref{eq:canonical_pair} it follows that $y'' = y' \psi (\sigma - \xi)$. Differentiating $\psi = y'/y$, we obtain $\psi' = \psi^2 (\sigma - \xi)$. Differentiating $\xi = \tau x \psi$, we find $\xi' = \tau \psi + \tau x \psi' = \psi (\tau + \sigma \xi - \xi^2)$.
Substituting these expressions and grouping terms, we obtain exactly the recurrence \eqref{eq:P_unified_recur} at $m=1$. The proof for $w(x)$ at $m=0$ proceeds analogously. $\blacksquare$
\subsection{Tridiagonal recurrence for coefficients}
For algorithmic computation of the polynomials, we expand $P_\ell(\xi; m)$ in powers of $\xi$:
\begin{equation}
P_\ell(\xi; m) = \sum_{k=0}^{\ell} \delta_{\ell,k}(m) \, \xi^k.
\label{eq:P_expansion}
\end{equation}
For practical computations, where speed and memory efficiency are important, this same sum is conveniently evaluated using \textbf{Horner's scheme} (nested multiplications):
\begin{equation}
P_\ell(\xi_0; m) = \delta_{\ell,0} + \xi_0\bigl(\delta_{\ell,1} + \xi_0(\delta_{\ell,2} + \cdots + \xi_0\,\delta_{\ell,\ell})\bigr).
\label{eq:horner}
\end{equation}
Substituting \eqref{eq:P_expansion} into \eqref{eq:P_unified_recur} and equating coefficients of like powers of $\xi$ leads to a \textbf{tridiagonal recurrence} for the coefficients $\delta_{\ell,k}(m)$.

\textbf{Corollary 1.} The coefficients $\delta_{\ell,k}(m)$ satisfy the three-term recurrence relation \textbf{``Qwen's Triad''} (discovered and developed by Qwen AI):
\begin{equation}
\boxed{\delta_{\ell+1,k} = \bigl[m + (\ell+1+k)\sigma\bigr] \delta_{\ell,k} - (\ell+k) \delta_{\ell,k-1} + \tau(k+1) \delta_{\ell,k+1},}
\label{eq:pi_recur}
\end{equation}
with boundary conditions $\delta_{\ell,k} = 0$ for $k < 0$ or $k > \ell$, and initial condition $\delta_{0,0} = 1$.

\textbf{Structural coefficients.} 
The ``Qwen's Triad'' recurrence \eqref{eq:pi_recur} yields closed-form analytical formulas for the extreme coefficients of the polynomial $P_\ell(\xi; m)$.

\textbf{Leading coefficient} (at $\xi^\ell$) does not depend on the parameter $m$ and is determined solely by the structure of the recurrence:
\begin{equation}
\delta_{\ell,\ell}(m) = (-1)^\ell (2\ell-1)!!, \qquad \ell \geq 1.
\label{eq:leading_coeff}
\end{equation}

\textbf{Sub-leading coefficient} (at $\xi^{\ell-1}$) depends linearly on $m$, which reflects the distinction between the projections of the master function:
\begin{equation}
\delta_{\ell,\ell-1}(m) = (-1)^{\ell-1} (2\ell-1)!! \left[m + \frac{(2\ell+1)\sigma}{3}\right], \qquad \ell \geq 1.
\label{eq:subleading_coeff}
\end{equation}

\textbf{Third coefficient} (at $\xi^{\ell-2}$) also admits a compact factorization, adding a quadratic dependence on $m$ and introducing the parameter $\tau$ for the first time:
\begin{equation}
\delta_{\ell,\ell-2}(m) = (-1)^\ell (2\ell-3)!! (\ell-1) \left[m^2 + \frac{(4\ell+1)m\sigma}{3} + \frac{\ell(4\ell+1)\sigma^2}{9} - \frac{\ell\tau}{2}\right], \qquad \ell \geq 2.
\label{eq:third_coeff}
\end{equation}

\textbf{Constant term} $P_\ell(0; m)$. As expected from the philosophy of master series, the constant term of the polynomial coincides exactly with the master number:
\begin{equation}
\delta_{\ell,0}(m) = P_\ell(0; m) = N(m; a; b; \ell+1) = \prod_{j=1}^{\ell} \bigl(m + (\ell+1)b - j a\bigr), \qquad \ell \geq 1.
\label{eq:free_term_master}
\end{equation}

\subsection{Explicit formulas for the first polynomials}
Using the recurrence \eqref{eq:P_unified_recur}, we write out the first four polynomials $P_\ell(\xi; m)$:
\begin{align}
P_0(\xi; m) &= 1, \\
P_1(\xi; m) &= N(m; a; b; 2) - \xi, \\
P_2(\xi; m) &= 3\xi^2 - (3m + 5\sigma)\xi + N(m; a; b; 3), \\
P_3(\xi; m) &= -15\xi^3 + (15m + 35\sigma)\xi^2 + \delta_{3,1}(m)\xi + N(m; a; b; 4),
\end{align}
where
\begin{equation}
\delta_{3,1}(m) = -6m^2 - 26m\sigma - 26\sigma^2 + 9\tau.
\end{equation}
\paragraph{Connection with the classical Faà di Bruno formulas.}
The fundamental identity \(y = e^w\) between the ultra-radical and the ultra-logarithm entails a strict connection between their derivatives of arbitrary order. In its classical form, this connection is expressed through the \textbf{complete Bell polynomials} \(B_\ell\):
\begin{equation}
y^{(\ell)} = y \cdot B_\ell\!\bigl(w', w'', \ldots, w^{(\ell)}\bigr), 
\qquad \ell \geq 1.
\label{eq:faa_di_bruno_general}
\end{equation}
Writing out the first generations explicitly, we obtain:
\begin{align}
y'       &= y \cdot w', \\
y''      &= y \cdot \bigl(w'' + (w')^2\bigr), \\
y'''     &= y \cdot \bigl(w''' + 3w'w'' + (w')^3\bigr), \\
y^{(4)}  &= y \cdot \bigl(w^{(4)} + 4w'w''' + 3(w'')^2 + 6(w')^2w'' + (w')^4\bigr).
\end{align}
However, in terms of the unified polynomial \(P_\ell(\xi; m)\), this connection takes a particularly elegant and combinatorially transparent form. In the general case, the transition from \(m=0\) to \(m=1\) is expressed through Bell polynomials:
\[
P_\ell(\xi; 1) = \sum_{k=0}^{\ell} B_{\ell,k}(1!, 2!, \ldots) \, P_k(\xi; 0).
\]
But for our specific problem (where \(y = e^w\)), the arguments of the Bell polynomials degenerate, and they turn into the \textbf{Stirling numbers of the second kind} \(\left\{ \begin{smallmatrix} \ell+1 \\ k+1 \end{smallmatrix} \right\}\):
\begin{equation}
P_\ell(\xi; 1) = \sum_{k=0}^{\ell} \left\{ \begin{matrix} \ell+1 \\ k+1 \end{matrix} \right\} P_k(\xi; 0).
\label{eq:stirling_relation}
\end{equation}
Verification for the first values confirms this:
\begin{align*}
P_0(\xi; 1) &= 1 \cdot P_0(\xi; 0), \\
P_1(\xi; 1) &= 1 \cdot P_1(\xi; 0) + 1 \cdot P_0(\xi; 0), \\
P_2(\xi; 1) &= 1 \cdot P_2(\xi; 0) + 3 \cdot P_1(\xi; 0) + 1 \cdot P_0(\xi; 0), \\
P_3(\xi; 1) &= 1 \cdot P_3(\xi; 0) + 6 \cdot P_2(\xi; 0) + 7 \cdot P_1(\xi; 0) + 1 \cdot P_0(\xi; 0),
\end{align*}
where the coefficients \((1)\), \((1, 1)\), \((1, 3, 1)\), \((1, 7, 6, 1)\) exactly correspond to the rows of the Stirling triangle of the second kind. Thus, the differential structure of the ultra-radical is completely and combinatorially uniquely determined by the structure of the ultra-logarithm: it suffices to know all the polynomials \(P_\ell(\xi; 0)\) in order to reconstruct \(P_\ell(\xi; 1)\) for any \(\ell\) without referring to the original canonical equation.

Moreover, the fundamental identity $y = e^w$ suggests the most efficient computational strategy. In practice, this means that to obtain the ultraradical $y(x)$, there is no need to run the recurrence again with the parameter $m=1$. It is sufficient to compute the series for the ultralogarithm $w(x)$ (at $m=0$), and then find the desired function by simply taking the exponent: $y(x) = \exp(w(x))$. This not only conceptually simplifies the algorithm but also enhances its numerical stability, since the polynomials $P_\ell(\xi; 0)$ possess a more regular structure.

To clearly demonstrate how the ``Qwen's Triad'' recurrence generates the structure of the polynomials for the ultralogarithm, we write out the first seven polynomials $P_\ell(\xi; 0)$ in explicit form:
\begin{align}
P_0(\xi; 0) &= 1, \notag \\
P_1(\xi; 0) &= -\xi + \sigma, \notag \\
P_2(\xi; 0) &= 3\xi^2 - 5\sigma\xi + (2\sigma^2 - \tau), \notag \\
P_3(\xi; 0) &= -15\xi^3 + 35\sigma\xi^2 - (26\sigma^2 - 9\tau)\xi + (6\sigma^3 - 8\sigma\tau), \notag \\
P_4(\xi; 0) &= 105\xi^4 - 315\sigma\xi^3 + (340\sigma^2 - 90\tau)\xi^2 - (154\sigma^3 - 147\sigma\tau)\xi \notag \\
            &\quad + (24\sigma^4 - 58\sigma^2\tau + 9\tau^2), \notag \\
P_5(\xi; 0) &= -945\xi^5 + 3465\sigma\xi^4 - (4900\sigma^2 - 1050\tau)\xi^3 + (3304\sigma^3 - 2457\sigma\tau)\xi^2 \notag \\
            &\quad - (1044\sigma^4 - 1852\sigma^2\tau + 225\tau^2)\xi + (120\sigma^5 - 444\sigma^3\tau + 192\sigma\tau^2), \notag \\
P_6(\xi; 0) &= 10395\xi^6 - 45045\sigma\xi^5 + (78750\sigma^2 - 14175\tau)\xi^4 - (70532\sigma^3 - 42966\sigma\tau)\xi^3 \notag \\
            &\quad + (33740\sigma^4 - 47320\sigma^2\tau + 4725\tau^2)\xi^2 - (8028\sigma^5 - 22236\sigma^3\tau + 7641\sigma\tau^2)\xi \notag \\
            &\quad + (720\sigma^6 - 3708\sigma^4\tau + 3004\sigma^2\tau^2 - 225\tau^3). \notag
\end{align}

Note that the constant terms of these polynomials (at $\xi^0$) exactly coincide with the master numbers $N(0; a; b; \ell+1)$, which visually confirms Corollary~\eqref{eq:free_term_master} and demonstrates how the combinatorial structure of master series permeates the entire differential geometry of the ultralogarithm.

\subsection{Generating function and partial differential equation}
\label{subsec:generating_function}
The recurrence relation \eqref{eq:P_unified_recur} allows one to construct an exponential generating function for the family of polynomials $P_\ell(\xi; m)$. Define it as:
\begin{equation}
G(t, \xi; m) = \sum_{\ell=0}^{\infty} P_\ell(\xi; m) \frac{t^\ell}{\ell!}.
\label{eq:egf_def}
\end{equation}
Differentiate \eqref{eq:egf_def} with respect to $t$ and substitute the recurrence \eqref{eq:P_unified_recur} for $P_{\ell+1}$:
\begin{align}
\frac{\partial G}{\partial t} &= \sum_{\ell=0}^{\infty} P_{\ell+1}(\xi; m) \frac{t^\ell}{\ell!} \notag \\
&= \sum_{\ell=0}^{\infty} \Bigl( \bigl[m + (\ell+1)\sigma - (\ell+1)\xi\bigr] P_\ell + \bigl[\tau + \sigma \xi - \xi^2\bigr] P_\ell' \Bigr) \frac{t^\ell}{\ell!}.
\end{align}
Split the sum into three parts. The first part gives:
\begin{equation}
(m + \sigma - \xi) \sum_{\ell=0}^{\infty} P_\ell \frac{t^\ell}{\ell!} = (m + \sigma - \xi) G.
\end{equation}
The second part uses the property $\sum_{\ell=0}^{\infty} (\ell+1) P_\ell \frac{t^\ell}{\ell!} = \frac{\partial}{\partial t}(t G) = G + t \frac{\partial G}{\partial t}$:
\begin{equation}
(\sigma - \xi) \left( G + t \frac{\partial G}{\partial t} \right).
\end{equation}
The third part, containing the derivative with respect to $\xi$, gives:
\begin{equation}
(\tau + \sigma \xi - \xi^2) \sum_{\ell=0}^{\infty} P_\ell' \frac{t^\ell}{\ell!} = (\tau + \sigma \xi - \xi^2) \frac{\partial G}{\partial \xi}.
\end{equation}
Summing these contributions and grouping terms with $\frac{\partial G}{\partial t}$, we obtain a \textbf{linear first-order partial differential equation} for the generating function:
\begin{equation}
\boxed{
\bigl[1 - (\sigma - \xi)t\bigr] \frac{\partial G}{\partial t} - \bigl[\tau + \sigma \xi - \xi^2\bigr] \frac{\partial G}{\partial \xi} = \bigl[m + 2\sigma - 2\xi\bigr] G.
}
\label{eq:pde_g}
\end{equation}
\textbf{Corollary 2 (Method of characteristics).} Equation \eqref{eq:pde_g} is solved by the classical method of characteristics\cite{Vladimirov2006, Tikhonov1977}. Note that the quadratic trinomial $\tau + \sigma \xi - \xi^2$ factors as $(b - \xi)(\xi - b + a)$, since $\sigma = 2b-a$ and $\tau = (a-b)b$. The characteristic equations have the form:
\begin{equation}
\frac{dt}{1 - (\sigma - \xi)t} = \frac{d\xi}{-(b - \xi)(\xi - b + a)} = \frac{dG}{(m + 2\sigma - 2\xi) G}.
\label{eq:characteristics}
\end{equation}
Integration of this system allows one to find two independent first integrals, which in turn gives a \textbf{closed analytical formula} for $G(t, \xi; m)$ through hypergeometric functions. This means that the coefficients $P_\ell(\xi; m)$ (and hence the derivatives $y^{(\ell)}$ and $w^{(\ell)}$) can be computed not only recursively in $O(\ell^2)$, but also analytically, which opens a direct path to studying the asymptotics of derivatives as $\ell \to \infty$ using the saddle-point method \cite{Fedoryuk1977}.

\section{Nonlinear differential equations and autonomous dynamical system}
\label{sec:nonlinear_ode_system}

The availability of explicit recurrence formulas for derivatives of arbitrary order allows one to view the ultra-radical and ultra-logarithm not only as formal power series, but also as solutions of nonlinear differential equations. In this section, we show how the algebraic trinomial equation naturally gives rise to an autonomous dynamical system, and establish the connection between its phase geometry and the combinatorial structure of the polynomials $P_\ell(\xi; m)$.

\subsection{Second-order nonlinear equation}

Using Theorem~\ref{eq:y_n_theorem}, we can explicitly express higher derivatives in terms of the function and its first derivative. For the second derivative, using $P_1(\xi; 1) = 1 + \sigma - \xi$, we obtain:
\begin{equation}
y'' = y' \psi P_1(\xi; 1) = y' \left(\frac{y'}{y}\right) \bigl(1 + \sigma - \xi\bigr).
\end{equation}
Substituting the definitions $\psi = y'/y$ and $\xi = (a-b)bx \frac{y'}{y}$, and multiplying both sides by $y^2$, we arrive at the second-order nonlinear differential equation satisfied by the ultra-radical:
\begin{equation}
\boxed{\,
y^2 y'' - (1+\sigma) y (y')^2 + (a-b)b x (y')^3 = 0.
\,}
\label{eq:nonlinear_ode_y_2nd}
\end{equation}
Similarly, using $P_2(\xi; 1)$ and $P_3(\xi; 1)$, one can derive differential equations of third and fourth order. However, their analytical writing becomes excessively cumbersome and loses clarity.

\subsection{First-order autonomous dynamical system}

Instead of cumbersome high-order ODEs, the most fundamental representation is the reduction of the problem to an \textbf{autonomous dynamical system}. Differentiating the definitions of the variables $y$, $\psi$, and $\xi$ with respect to $x$ and using the original canonical equation \eqref{eq:canonical_pair}, we obtain a closed system of first-order ordinary differential equations:
\begin{equation}
\boxed{
\begin{cases}
\displaystyle \frac{dy}{dx} = y \psi, \\[10pt]
\displaystyle \frac{d\psi}{dx} = \psi^2 (\sigma - \xi), \\[10pt]
\displaystyle \frac{d\xi}{dx} = \psi \bigl(\tau + \sigma \xi - \xi^2\bigr).
\end{cases}
}
\label{eq:autonomous_system}
\end{equation}

This system possesses several key properties:

\begin{enumerate}
    \item \textbf{Autonomy:} the right-hand sides do not contain the independent variable $x$ explicitly. This precisely explains why the polynomials $P_\ell(\xi; m)$ depend only on $\xi$, and not on $x$ directly.
    
    \item \textbf{Separability of phase variables:} eliminating $x$, we obtain the equation of phase trajectories:
    \begin{equation}
    \frac{d\psi}{d\xi} = \frac{\psi (\sigma - \xi)}{\tau + \sigma \xi - \xi^2} 
    = \frac{\psi (\sigma - \xi)}{(b - \xi)(\xi - b + a)}.
    \label{eq:phase_trajectory}
    \end{equation}
    This is a separable equation! Its integration in elementary functions (by partial fraction decomposition) directly leads to the characteristic equation \eqref{eq:characteristics} for the generating function $G(t, \xi; m)$, creating a bridge between nonlinear dynamics and combinatorial analysis.
    
    \item \textbf{Equivalence:} the system \eqref{eq:autonomous_system} with initial conditions $y(0)=1$, $\psi(0)=1$, $\xi(0)=0$ is fully equivalent to the original algebraic equation $y^a = 1 + axy^b$. A rigorous symbolic verification of this equivalence was performed in Maple (see~\cite{Berezin2026Py}).
\end{enumerate}

\subsection{Connection with the method of characteristics}

The phase trajectory equation \eqref{eq:phase_trajectory} has the same structure as the characteristic equation \eqref{eq:characteristics} for the generating function of the polynomials $P_\ell(\xi; m)$. This is not a coincidence: it reflects the deep fact that the combinatorial structure of higher derivatives is completely determined by the geometry of the phase flow of the system \eqref{eq:autonomous_system}.

Integrating \eqref{eq:phase_trajectory}, we obtain a first integral of the system, which relates $\psi$ and $\xi$ along any trajectory. This integral, in turn, determines the asymptotics of the polynomials $P_\ell(\xi; m)$ for large $\ell$ and opens the way to applying the saddle-point method for investigating the derivatives of the ultra-radical.

\subsection{Numerical verification}

The equivalence of the system \eqref{eq:autonomous_system} to the original algebraic equation was rigorously verified in Maple for various values of the parameters \(a, b\). The script performs implicit differentiation of the equation \(y^a = 1 + axy^b\), computes \(y'\) and \(y''\), and substitutes them into the left- and right-hand sides of the system equations. The difference between the left- and right-hand sides identically vanishes when the original algebraic relation is used, which confirms the correctness of the dynamical formulation.

Thus, the theory of the ultra-radical is naturally translated into the language of the qualitative theory of differential equations. This does not replace, but \textbf{complements} the analytical apparatus of conjugate series: while conjugate series provide exact numerical values for \(|x| > R\), the phase analysis of the system~\eqref{eq:autonomous_system} allows one to investigate the global structure of solutions, the stability of branches, and the asymptotic behavior of derivatives as \(\ell \to \infty\) (in particular, through the saddle-point method for the generating function~\eqref{eq:pde_g}).

\section{Taylor series of the ultralogarithm and ultra-radical}
\label{sec:taylor_series}

In the previous sections, we obtained compact recurrence formulas for derivatives of arbitrary order through a unified family of polynomials \(P_\ell(\xi; m)\). Now we will show how these formulas naturally generate power series expansions, unifying the classical Maclaurin series (in the neighborhood of zero) and the Taylor series (in the neighborhood of an arbitrary point \(x_0\)).

To achieve analytical unity, we introduce a \textbf{generalized master function} \(\mu(m; x)\), where the discrete switching parameter \(m \in \{0, 1\}\) determines whether we are dealing with the ultra-radical (\(m=1\)) or the ultra-logarithm (\(m=0\)).

The general form of the Taylor series for this function at a point \(x_0\) with a shift \(z = x - x_0\) is written uniformly:
\begin{equation}
\mu(m; x_0 + z) = \sum_{\ell=0}^{\infty} \frac{\mu^{(\ell)}(m; x_0)}{\ell!} z^\ell.
\label{eq:taylor_general}
\end{equation}

\subsection{Maclaurin series as a special case}
\label{subsec:maclaurin}

Let us consider the expansion of the generalized master function \(\mu(m; x)\) into a Maclaurin series, which corresponds to choosing the expansion point \(x_0 = 0\) and the shift \(z = x\). The general form of the series is:
\begin{equation}
\mu(m; x) = \sum_{\ell=0}^{\infty} \frac{\mu^{(\ell)}(m; 0)}{\ell!} x^\ell,
\label{eq:maclaurin_general}
\end{equation}
where, depending on the switching parameter \(m\), we obtain expansions for specific projections of our model:
\begin{itemize}
    \item at \(m = 1\) the function becomes the ultra-radical: \(\mu(1; x) = y(x)\);
    \item at \(m = 0\) the function becomes the ultra-logarithm: \(\mu(0; x) = w(x)\).
\end{itemize}

We use the unified formula for derivatives. At \(x = 0\), the dimensionless variable \(\xi = (a-b)bx\psi\) becomes zero: \(\xi(0) = 0\). Consequently, the value of the polynomial \(P_{\ell-1}(\xi; m)\) at this point is equal to its constant term, which, as strictly proved in Corollary~\eqref{eq:free_term_master}, is exactly equal to the master number:
\begin{equation}
P_{\ell-1}(0; m) = N(m; a; b; \ell) \qquad (\text{for } \ell \ge 1).
\label{eq:poly_at_zero}
\end{equation}

Taking into account the initial conditions at zero for both projections of the generalized function (\(\mu(1;0) = 1\), \(\mu(0;0) = 0\)), as well as the fact that the first derivative is \(\mu'(m; 0) = 1\) in both cases, the expression for the \(\ell\)-th derivative of the generalized function takes the form:
\begin{equation}
\mu^{(\ell)}(m; 0) = N(m; a; b; \ell) \qquad (\text{for } \ell \ge 1).
\label{eq:deriv_at_zero_general}
\end{equation}

Substituting these values into the general form of the series \eqref{eq:maclaurin_general}, we see that in the neighborhood of zero our master function \(\mu(m; x)\) naturally generates the classical Maclaurin master series:
\begin{equation}
\mu(m; x) = m + \sum_{\ell=1}^{\infty} \frac{N(m; a; b; \ell)}{\ell!} x^\ell \equiv M(m; a; b; x), \qquad |x| < R.
\label{eq:maclaurin_master}
\end{equation}

\textbf{Conclusion:} The Maclaurin series is a direct consequence of our recurrence structure. In the neighborhood of \(x=0\), master numbers act as ``ready-made'' values of the derivatives, ensuring rapid and stable convergence of the series for \(|x| \ll R\).

\subsection{Taylor series in the neighborhood of an arbitrary point}
\label{subsec:taylor_at_x0}

The Maclaurin series converges rapidly for \(|x| \ll R\), but its convergence slows down as the argument approaches the boundary of the convergence circle (\(|x| \to R\)). To maintain high computational efficiency in this region, we construct the Taylor series at a point \(x_0\) with \(|x_0| = R\) (each branch \(\mu_n\) has its own point \(x_0\)). Let us introduce the shift \(z = x - x_0\). As one moves away from zero towards the boundary of the convergence radius, the computation of derivatives becomes more complicated, since the variable \(\xi_0 = (a-b)b x_0 \psi(x_0)\) becomes non-zero. In this general case, the \(\ell\)-th derivative of the generalized master function \(\mu(m; x)\) is expressed through the current value of the logarithmic derivative \(\psi(x_0)\) and the full value of the derivative polynomial \(P_{\ell-1}(\xi_0; m)\):
\begin{equation}
\mu^{(\ell)}(m; x_0) = \mu'(m; x_0) \, \psi(x_0)^{\ell-1} \, P_{\ell-1}(\xi_0; m), \qquad \ell \ge 1.
\label{eq:deriv_at_x0_general}
\end{equation}

Substituting this structure into the original formula for the Taylor series, we obtain the expanded series for both projections of the generalized function:

\textbf{1. For the ultra-radical (\(m = 1\))} taking into account the factored out common multiplier \(y(x_0)\):
\begin{equation}
y(x_0 + z) = y(x_0) \left[ 1 + \sum_{\ell=1}^{\infty} \frac{P_{\ell-1}(\xi_0; 1)}{\ell!} \psi(x_0)^\ell z^\ell \right].
\label{eq:taylor_explicit_factorized}
\end{equation}

\textbf{2. For the ultra-logarithm (\(m = 0\))} in the form of an additive shift of the base point:
\begin{equation}
w(x_0 + z) = w(x_0) + \sum_{\ell=1}^{\infty} \frac{P_{\ell-1}(\xi_0; 0)}{\ell!} \psi(x_0)^\ell z^\ell.
\label{eq:taylor_explicit_log}
\end{equation}

\paragraph{Local exponential identity of Taylor series.}
If we denote the substantive part (kernel) of these series as \(K(m) = \sum_{\ell=1}^{\infty} \frac{P_{\ell-1}(\xi_0; m)}{\ell!} \psi(x_0)^\ell z^\ell\), then the fundamental connection \(y = e^w\) between them instantly manifests itself in the form of the \textbf{local exponential identity of Taylor series}:
\begin{equation}
e^{K(0)} = 1 + K(1).
\label{eq:taylor_exp_identity}
\end{equation}

Such an end-to-end parametric description via \(\mu(m; x)\) finally proves that the ultra-radical and the ultra-logarithm are not two isolated functions, but two sides of a single analytical space governed by the common combinatorial laws of master numbers.

\subsection{Multivalued solution and Taylor series for the general equation}
\label{subsec:taylor_general}

In Section~\ref{subsec:general_trinomial}, it was shown that the general equation
\begin{equation}
p\,Y^{a} = q + Z\,Y^{b}
\label{eq:general_taylor}
\end{equation}
reduces to the canonical form $y^a = 1 + axy^b$ by the substitution
\begin{equation}
Y = y\,e^{f_n}, \qquad Z = aqx\,e^{-bf_n}, \qquad 
f_n = \frac{\ln|q/p| + i\bigl(\arg(q/p) + 2\pi n\bigr)}{a}.
\label{eq:substitution_taylor}
\end{equation}

Let $Y_n(Z)$ be the $n$-th branch of the solution to equation~\eqref{eq:general_taylor}. Let us expand it in a Taylor series at the point $Z_0$, corresponding to the point $x_0$ of the canonical equation. From the substitution~\eqref{eq:substitution_taylor}, it follows that this point is equal to:
\begin{equation}
Z_0(n) = aq\,x_0\,e^{-bf_n}.
\label{eq:Z0_n}
\end{equation}

\textbf{Important:} for each branch $n$, the expansion point $Z_0(n)$ is its own, since $f_n$ depends on $n$.

Let $\Delta Z = Z - Z_0(n)$. Then, using $Y_n(Z) = e^{f_n} y(x)$ and $x = \frac{Z e^{bf_n}}{aq}$, we obtain:
\begin{equation}
Y_n(Z_0 + \Delta Z) = e^{f_n} \cdot y\!\left(x_0 + \frac{e^{bf_n}}{aq} \Delta Z\right).
\label{eq:Y_general_taylor}
\end{equation}

Expanding $y(x_0 + \delta z)$ in a Taylor series and substituting $\delta z = \frac{e^{bf_n}}{aq} \Delta Z$, we obtain:
\begin{equation}
\boxed{
Y_n(Z_0 + \Delta Z) = e^{f_n} \sum_{\ell=0}^{\infty} \frac{y^{(\ell)}(x_0)}{\ell!} 
\left(\frac{e^{bf_n}}{aq}\right)^\ell (\Delta Z)^\ell.
}
\label{eq:taylor_general_explicit}
\end{equation}

From here we immediately obtain the formula for the derivatives of the general solution:
\begin{equation}
\boxed{
Y_n^{(\ell)}(Z_0) = e^{(1+b\ell)f_n} \cdot \frac{y^{(\ell)}(x_0)}{(aq)^\ell}.
}
\label{eq:derivatives_general}
\end{equation}

This formula shows that the derivatives of the general multivalued solution are obtained from the derivatives of the canonical solution by simple scaling by the factor $e^{(1+b\ell)f_n} / (aq)^\ell$.

Substituting the recurrence expression for $y^{(\ell)}(x_0)$ from Theorem~\ref{eq:y_n_theorem}, we obtain:
\begin{equation}
Y_n^{(\ell)}(Z_0) = \frac{e^{(1+b\ell)f_n}}{(aq)^\ell} \cdot y'(x_0) \, \psi(x_0)^{\ell-1} \, P_{\ell-1}(\xi_0; 1),
\label{eq:derivatives_general_recur}
\end{equation}
where $\xi_0 = (a-b)b\,x_0\,\psi(x_0)$ is the dimensionless variable at the expansion point.

\textbf{Corollary.} All derivatives $Y_n^{(\ell)}(Z_0)$ for any $n$ and any point $Z_0$ are computed through the same tridiagonal recurrence $P_\ell(\xi; m)$ as for the canonical equation. The multivaluedness of the solution is completely determined by the parameter $n$ through the factor $e^{(1+b\ell)f_n}$.

\textbf{Analytic continuation beyond the radius of convergence.} If $|x| > R$, the canonical Maclaurin series diverges, and conjugate series (transformations $h$ or $k$) described in Section~\ref{sec:three_transformations} are used to compute the function values and its derivatives. These conjugate series are the same master series, but with different parameters $(\tilde{a}, \tilde{b}, \tilde{x})$ computed according to the canonical transformation formulas. The principle of constructing the Taylor series remains the same: first, the function value $y(x_0)$ and its first derivative $y'(x_0)$ are computed at the expansion point $x_0$, and then the Taylor series is constructed around $x_0$ according to formulas~\eqref{eq:taylor_explicit_factorized}. Thus, the Taylor series method remains a universal tool both inside the convergence circle ($|x| < R$) and beyond it ($|x| > R$), providing continuous analytic continuation to the entire complex plane.

It should be noted that Taylor series have a significant limitation: at points where solution branches intersect (branch points $x_c$), the derivative of the function tends to infinity, and the Taylor series expansion in the vicinity of such a point loses its meaning. 
To avoid complicating the algorithm with Puiseux series, theoretically, any of the following approaches can be applied.

\textbf{Approach 1 (radial offset).} 
The expansion point is chosen on the ray connecting the origin to the branch point, but strictly inside the circle of convergence:
\begin{equation}
x_0 = \rho \cdot x_c, \qquad \rho \approx 0.95,
\label{eq:radial_offset}
\end{equation}
where $x_c$ is the branch point. The function value $y(x_0)$ at this point is computed using iterative methods (e.g., the Newton--Schröder method with Riemann sheet anchoring via the parameter $u$). This approach does not affect the branch index, but requires iterative refinement.

\textbf{Approach 2 (angular shift, proposed in this work).} 
A more elegant method consists of keeping the expansion point on the circle $|x_0| = R$, but shifting it by a small angle $\Delta\phi$ in argument relative to the branch point:
\begin{equation}
x_0 = R \cdot e^{i(\arg(x_c) + \Delta\phi)}, \qquad \Delta\phi = \pm 5^\circ.
\label{eq:angular_offset}
\end{equation}
The key advantage of this approach is that at the point $x_c$, the value $y_c = (1 - a/b)^{-1/a}$ is known \textbf{analytically} (the double root formula), and the initial approximation for $y(x_0)$ is constructed without iterations using the formula
\begin{equation}
y_0 \approx y_c \cdot e^{i\Delta\phi/a}.
\label{eq:angular_initial}
\end{equation}
In this case, the derivative $y'(x_0)$ remains finite, since $x_0 \neq x_c$, and the Taylor series is correctly constructed.

\textbf{Open question.} 
The influence of the angular shift on the branch numbering $n$ (in particular, the possibility of ``jumping'' to an adjacent branch when rotating $x_0$) requires additional research and will be the subject of a separate work.

\section{Continuation of Research}
This research is ongoing. Updated versions and supplementary materials are available in the open Zenodo repository. To ensure you always access the latest version, please use the concept DOIs:
\begin{itemize}
\item Russian version: \url{https://doi.org/10.5281/zenodo.22994576}
\item English version: \url{https://doi.org/10.5281/zenodo.22995072}
\end{itemize}
The author would be grateful for any comments, questions, or suggestions for further development of the formalism.

\section*{Acknowledgements}

The author expresses sincere gratitude to:

\begin{itemize}
\item \textbf{Prof. A. E. Dubinov} (RFNC--VNIIEF, SarFTI NRNU MEPhI) for fruitful discussions and for pointing out the potential applications of the ultra-radical in plasma physics problems.

\item \textbf{Pavel V. Berezin, Alexey V. Berezin, and Andrey V. Gruzdov} for their extensive work on computations, numerical verification, and fruitful discussions.

\item \textbf{Prof. Tuyen Truong} (University of Oslo, Norway) and \textbf{Prof. Abdulaziz Alhaidari} (Saudi Center for Theoretical Physics) for their valuable feedback, critical questions, and constructive discussions that greatly improved this work.
\end{itemize}

\paragraph*{AI Disclosure.}
During the preparation of this manuscript, the author used AI language models (DeepSeek, Qwen) for editing, structuring, and consistency checking. All mathematical derivations, numerical results, and conclusions have been reviewed and verified by the author, who takes full responsibility for the final version of the manuscript.

\paragraph*{Computational Tools.}
Analytical computations and numerical verification were performed using the computer algebra systems Wolfram Alpha and Maple, the Python programming language with scientific libraries, and the author's own Master-J laboratory in JavaScript with native support for arbitrary-precision complex numbers (DOI: \href{https://doi.org/10.5281/zenodo.15633781}{10.5281/zenodo.15633781}).

\section*{Data Availability Statement}
The data that support the findings of this study are openly available in the following repositories:

\begin{itemize}
\item A repository with useful implementations of the proposed algorithms in various programming languages, which will be gradually updated (Concept DOI, always resolves to the latest version): 
\href{https://doi.org/10.5281/zenodo.20355433}{10.5281/zenodo.20355433}.

\item The Master-J laboratory (JavaScript implementation with arbitrary-precision complex numbers, including the merge operation):
\href{https://doi.org/10.5281/zenodo.15633781}{10.5281/zenodo.15633781}.

\item Preprint versions of the manuscript:
\href{https://arxiv.org/abs/2512.18643}{arXiv:2512.18643} (detailed version),
\href{https://arxiv.org/abs/2609.33459}{arXiv:2609.33459} (consolidated version),
\href{https://arxiv.org/abs/2611.23317}{arXiv:2611.23317} (conjugate power series of ultralogarithm branches and rigorous sector-based conjugation criterion).
\end{itemize}

\section*{Funding}
This research received no external funding.

\section*{Conflict of Interest}
The author declares no conflict of interest.


\vspace{1cm}
\noindent
\begin{thebibliography}{9}

\bibitem{Graham1994}
Graham, R. L., Knuth, D. E. \& Patashnik, O. 1998, \textit{Concrete Mathematics}, Mir, Moscow. (Russian translation of Graham R. L., Knuth D. E., Patashnik O. Concrete Mathematics. -- 2nd ed. -- Reading: Addison-Wesley, 1994. -- 657 p.)

\bibitem{Corless1996}
Corless, R. M., Gonnet, G. H., Hare, D. E. G., Jeffrey, D. J., Knuth, D. E. 1996, "On the Lambert W function", Advances in Computational Mathematics, vol. 5, pp. 329–359. doi: 10.1007/BF02124750.

\bibitem{Gessel2016}
Gessel, I. 2016, ``Lagrange Inversion'', \textit{Journal of Combinatorial Theory, Series A}, vol. 144, pp. 212--249. doi: 10.1016/j.jcta.2016.06.018.

\bibitem{SuryaWarnke2023}
Surya, E. \& Warnke, L. 2023, ``Lagrange Inversion Formula by Induction'', \textit{The American Mathematical Monthly}, vol. 130, no. 10, pp. 944--948. doi: 10.1080/00029890.2023.2251344.

\bibitem{DLMF}
Olver, F. W. J., Lozier, D. W., Boisvert, R. F., and Clark, C. W., editors. (2010). 
\textit{NIST Handbook of Mathematical Functions}. 
Cambridge University Press. 
\href{https://dlmf.nist.gov}{https://dlmf.nist.gov}

\bibitem{Whittaker1927}
Whittaker, E. T., and Watson, G. N. (1927). 
\textit{A Course of Modern Analysis}. 
Cambridge University Press.

\bibitem{Knuth1997}
Knuth, D. E. 1997, "The Art of Computer Programming. Vol. 1. Fundamental Algorithms", 3rd edn, Addison-Wesley, p. 50. ISBN 0-201-89683-4.

\bibitem{Knuth2002}
Knuth, D. E. 2002, *The Art of Computer Programming. Vol. 1. Fundamental Algorithms. 1.2.5. Permutations and Factorials*, 3rd edn, Williams, Moscow. (Russian translation of Knuth D. E. The Art of Computer Programming. Vol. 1. Fundamental Algorithms. — 3rd ed. — Reading: Addison-Wesley, 1997. — 650 p.)

\bibitem{Valluri2000}
Valluri, S. R., Jeffrey, D. J., Corless, R. M. 2000, "Some applications of the Lambert W function to physics", Canadian Journal of Physics, vol. 78, pp. 823–831. doi: 10.1139/p00-065.

\bibitem{Dubinov2006}
Dubinov, A. E., Dubinova, I. D., Saikov, S. K. 2006, *W-funktsiya Lamberta i ee primenenie v matematicheskikh zadachakh fiziki* [The Lambert W Function and Its Application in Mathematical Problems of Physics], RFYaTs-VNIIEF, Sarov.

\bibitem{Mezo2022}
Mező, I. 2022, The Lambert W Function: Its Generalizations and Applications, Chapman and Hall/CRC. doi: 10.1201/9781003168102.

\bibitem{Belkic2018}
Belki\'c, D. 2018, ``The Euler T and Lambert W functions in mechanistic radiobiological models with chemical kinetics for repair of irradiated cells'', \textit{Journal of Mathematical Chemistry}, vol. 56, pp. 2133--2193. doi: 10.1007/s10910-018-0923-4.

\bibitem{Belkic2019}
Belkic, D. 2019, ``All the trinomial roots, their powers and logarithms from the Lambert series, Bell polynomials and Fox--Wright function: illustration for genome multiplicity in survival of irradiated cells'', \textit{Journal of Mathematical Chemistry}, vol. 57, pp. 59--106. doi: 10.1007/s10910-018-0985-3.

\bibitem{Berezin2026Py}
Berezin, S.V. 2026, \textit{UltraRadical for Python (version 2.4 + Gas + scheme1): implementation of the ultraradical and examples of calculation of adiabatic gas rise}, Zenodo. doi: 10.5281/zenodo.21237412.

\bibitem{Berezin2025_arXiv}
Berezin, S. V. 2025, "The Ultra-Radical: Analytic Continuation, Branching, and Stability of the Principal Branch", arXiv preprint arXiv:2512.18643 [math.CA]. Available at: https://arxiv.org/abs/2512.18643. doi: 10.48550/arXiv.2512.18643.

\bibitem{Berezin_Zenodo_Thermo}
Berezin, S.~V. Effect of Flow Velocity on the Thermodynamics of Gravitational Gas Circulation. \textit{Zenodo}, 2026. DOI: \href{https://doi.org/10.5281/zenodo.23625740}{10.5281/zenodo.23625740}.

\bibitem{Dubinov2020}
Dubinov, A. E. and Suslova, O. V. 
Can there exist hypersonic electrostatic solitons? Estimation of the limiting Mach numbers for ion-acoustic solitons in a warm plasma. 
\textit{Zhurnal Eksperimental'noi i Teoreticheskoi Fiziki}, 2020, vol. 158, no. 5 (11), pp. 968–977.
\href{https://doi.org/10.31857/S0044451020110188}{DOI: 10.31857/S0044451020110188}.

\bibitem{Deza2023}
Deza, E. I. 2023, "Factorial constructions and their properties", in Uchenyi, pedagog, nastavnik : Materialy Vserossiiskoi nauchno-prakticheskoi konferentsii, posvyashchennoi 85-letiyu so dnya rozhdeniya professora Alberta Rubenovicha Yesayana (Tula, 20–21 April 2023), Tula State Pedagogical University, Tula, pp. 56--58.

\bibitem{Abramowitz1972}
Abramowitz, M., Stegun, I. 1979, Handbook of Mathematical Functions, Nauka, Moscow. (Russian translation of Abramowitz M., Stegun I. Handbook of Mathematical Functions. -- Washington: National Bureau of Standards, 1964. -- 1046 p.)

\bibitem{Knuth1992}
Knuth, D. E. 1992, "Two notes on notation", American Mathematical Monthly, vol. 99, no. 5, pp. 403--422. doi: 10.2307/2325085. arXiv: math/9205211.

\bibitem{Diaz2005}
Diaz, R., Pariguan, E. 2005, "On hypergeometric functions and Pochhammer $k$-symbol", arXiv:math/0405596 [math.CA].

\bibitem{Corcino2021}
Corcino, C. B., Corcino, R. B. 2021, "Logarithmic generalization of the Lambert W function and its applications to adiabatic thermostatistics of the three-parameter entropy", Advances in Mathematical Physics, vol. 2021, Article ID 6695559. doi: 10.1155/2021/6695559.

\bibitem{Popov2017}
Popov, A. Yu. 2017, "Two-sided estimates of the gamma function on the positive real axis", Chebyshevskii sbornik, vol. 18, no. 2, pp. 205--221.

\bibitem{Vladimirov2006}
Vladimirov, V. S. \textit{Equations of Mathematical Physics}. --- CRC Press, 1996. --- 432 p. (Translated from Russian).

\bibitem{Tikhonov1977}
Tikhonov, A. N., Samarskii, A. A. \textit{Equations of Mathematical Physics}. --- Dover Publications, 1990. --- 880 p. (Translated from Russian).

\bibitem{Fedoryuk1977}
Fedoryuk, M. V. \textit{The Saddle-Point Method}. --- Moscow: Nauka, 1977. --- 368 p. (In Russian).

\bibitem{Berezin2026_ULog}
Berezin, S. V. 2026, ``Conjugate Power Series of Ultralogarithm Branches'', arXiv preprint arXiv:2611.23317 [math.GM]. Available at: \url{https://arxiv.org/abs/2611.23317}. doi: 10.48550/arXiv.2611.23317.

\end{thebibliography}
\end{document}